%% file: Elastic_Splines_II.tex
\magnification\magstep1
\overfullrule=0pt
\voffset=-0.3 in
\input amstex
\documentstyle{amsppt}
\input SecThEq

\input epsf.tex
\def\gk{\kappa}

\def\x{\xi}
\def\gep{\varepsilon}
\def\ga{\alpha}
\def\gb{\beta}
\def\gga{\gamma}
\def\gG{\Gamma}
\def\gD{\Delta}
\def\gs{\sigma}

\def\gt{\tau}

\def\gO{\Omega}
\def\Sopt{S_{\text{opt}}}
\def\tbar{\overline t}
\def\psibar{\overline \psi}

\def\BbR{\Bbb R}
\def\BbC{\Bbb C}
\def\BbZ{\Bbb Z}

\def\splice{\sqcup}

\def\norm#1{{\Vert #1 \Vert}}
\def\set#1{\{#1\}}
\def\abs#1{\left\vert {#1} \right\vert}
\def\pr#1{\left( {#1} \right)}
\def\wt{\widetilde}
\def\wh{\widehat}
\def\dd{\,{\roman d}}
\def\bs{\backslash}

\def\figsixteen{1}
\def\figseventeen{2}
\def\figtwelve{3}
\def\sunique{3.1}
\def\SectionRES{4}
\def\consistent{4.1}
\def\horseman{4.3}
\def\Psidefine{4.2}
\def\refinedcondCtwo{4.4}
\def\nothingyet{4.6}
\def\figfifteen{5}
\def\SectionQ{5}
\def\detq{5.3}
\def\cortstar{5.4}
\def\cortbar{5.5}
\def\SectionQunique{6}
\def\nablaEQ{7.1}

\def\AmProdi{1}
\def\BirkBoor{2}
\def\BJone{3}
\def\Brunnett{4}
\def\FisherJerome{5}
\def\chwuho{6}
\def\JeromeOne{7}
\def\JeromeTwo{8}
\def\JJ{9}
\def\LeeForsyth{10}
\def\LinnerJerome{11}

\hsize=6.5truein
\vsize=9truein

\topmatter

\title     Elastic Splines II: unicity of optimal s-curves and $G^2$ regularity of splines
\endtitle
    \author Albert Borb\'ely  \& Michael J. Johnson   \endauthor

    \address Department of Mathematics, Faculty of
    Science, Kuwait University, P.O. Box 5969, Safat 13060, Kuwait
    \endaddress

\email borbely.albert\@gmail.com, yohnson1963\@hotmail.com  \endemail
\thanks This work was supported and funded by Kuwait University, Research Project No.\ SM 01/14 
\endthanks

    \keywords  spline, nonlinear spline, elastica, bending energy, curve fitting, interpolation
    \endkeywords
    \subjclass   41A15, 65D17, 41A05    \endsubjclass
    \abstract
Given points $P_1,P_2,\ldots,P_m$ in the complex plane, we are concerned with the problem of finding
an interpolating curve with minimal bending energy (i.e., an optimal interpolating curve). It was shown
previously that existence is assured if one requires that the pieces of the interpolating curve be s-curves.
In the present article we also impose the restriction that these s-curves have chord angles not exceeding
$\pi/2$ in magnitude. 
With this setup, we have identified a sufficient condition for the $G^2$ regularity of optimal interpolating
curves. This sufficient condition relates to the stencil angles $\set{\psi_j}$, where $\psi_j$ is defined
as the angular
change in direction from segment $[P_{j-1},P_j]$ to segment $[P_j,P_{j+1}]$. A distinguished angle
$\Psi$ ($\approx 37^\circ$) is identified, and we show that if the stencil angles satisfy $\abs{\psi_j}<\Psi$, then
optimal interpolating curves are globally $G^2$.

As with the previous article, most of our effort is concerned with the geometric Hermite interpolation problem 
of finding an optimal s-curve which
connects $P_1$ to $P_2$ with prescribed chord angles $(\ga,\gb)$. Whereas existence 
was previously shown, and
sometimes uniqueness, the present article begins by establishing uniqueness when $\abs\ga,\abs\gb\leq\pi/2$ and
$\abs{\ga-\gb}<\pi$.
   \endabstract
\leftheadtext{Elastic Splines II}
\rightheadtext{unicity of optimal s-curves and $G^2$ regularity of splines}
\endtopmatter
\document
%
\Section{Introduction}
%
Given points $P_1,P_2,\ldots,P_m$ in the complex plane $\BbC$ with $P_j\neq P_{j+1} $, we are concerned with
the problem of finding a {\it fair} curve which interpolates the given points. The present contribution
is a continuation of \cite{\BJone} and so we adopt much of the notation used there. In particular, an
{\bf interpolating curve} is an absolutely-continuously differentiable function $F:[a,b]\to\BbC$,
with $F'$ non-vanishing, for
which there exist times $a=t_1<t_2<\cdots<t_m=b$ such that $F(t_j)=P_j$. We treat $F$ as a curve consisting
of $m-1$ pieces; the $j$-th piece of $F$, denoted $F_{[t_j,t_{j+1}]}$, runs from $P_j$ to $P_{j+1}$.
It is known (see \cite{\BirkBoor})
that there does not exist an interpolating curve with minimal bending energy, except in the trivial case
when the interpolation points lie sequentially along a line. In \cite{\BJone}, it was shown that existence
is assured if one imposes the additional condition that each piece of the interpolating curve 
be an s-curve.  Here, an {\bf s-curve} is a curve which first turns
monotonically at most $180^\circ$ in one
direction (either counter-clockwise or clockwise) and then turns monotonically at most $180^\circ$ in
the opposite direction. Incidentally, a {\bf c-curve} is an s-curve which turns in only one direction, 
and a {\bf u-turn} is a c-curve which turns a full $180^\circ$.
Associated with an s-curve $f:[a,b]\to \BbC$ (see Fig.\ \figsixteen) are its {\bf breadth} $L=\abs{f(b)-f(a)}$ 
and {\bf chord angles} $(\ga,\gb)$, defined by
$$
\ga=\arg \frac{f'(a)}{f(b)-f(a)},\quad \gb=\arg \frac{f'(b)}{f(b)-f(a)},
$$
where $\arg$ is defined with the usual range $(-\pi,\pi]$. \newline
\epsfysize=0.5truecm \epsffile{space.eps}\epsfysize=3.2truecm \epsffile{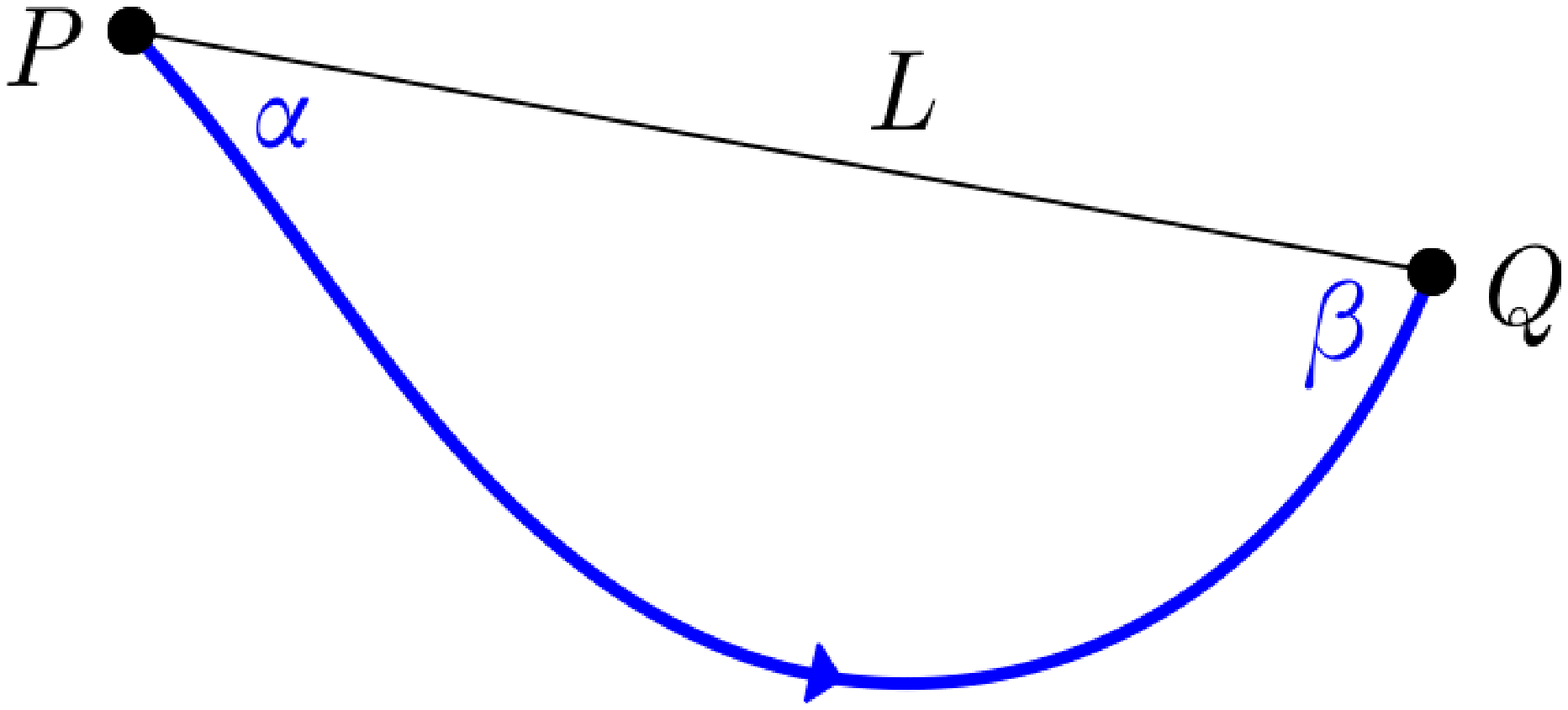}\hskip1cm
\epsfysize=3.2truecm \epsffile{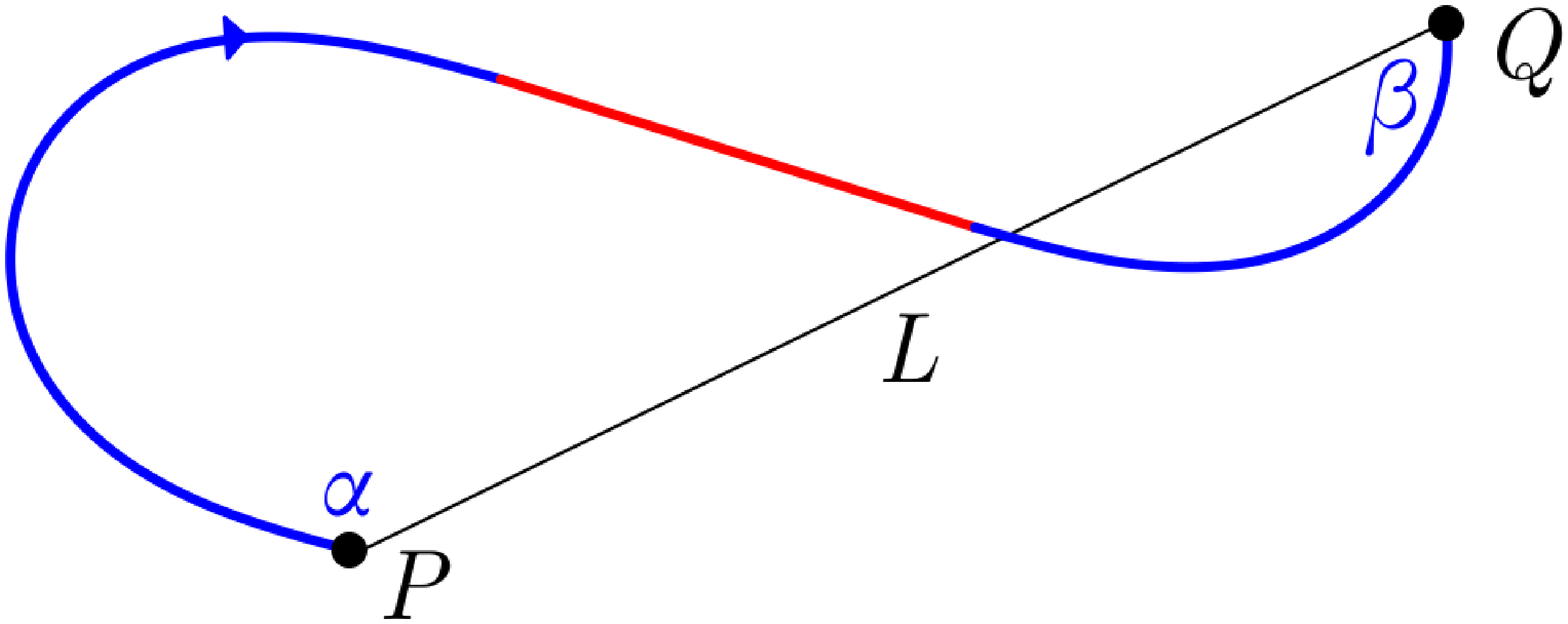}\newline
{\bf Fig.\ \Flabel\figsixteen\ \ (a)} optimal s-curve of Form 1 \ \ \ \ \ \ \ \ \ \ \ \ 
{\bf (b)} optimal s-curve of Form 2\newline
Note that although the chord angles are signed, our figures only indicate their magnitudes.
The chord angles $(\ga,\gb)$ of an s-curve necessarily satisfy
$$
|\ga|,|\gb|<\pi\quad\text{ and }\quad |\ga-\gb|\leq \pi.
\tag\Elabel\scurveineq$$

Defining $$\Cal A(P_1,P_2,\ldots,P_m)$$ to be the set of all interpolating curves
whose pieces are s-curves, the main result of \cite{\BJone} is that $\Cal A(P_1,P_2,\ldots,P_m)$ contains
a curve (called an {\bf elastic spline})
with minimal bending energy. Most of the effort in \cite{\BJone} is devoted to proving the existence of optimal s-curves.
Specifically, it is shown that given distinct points $P,Q$ and angles $(\ga,\gb)$ satisfying (\scurveineq), the set of all
s-curves from $P$ to $Q$ with chord angles $(\ga,\gb)$ contains a curve with minimal bending energy. 
Denoting the bending energy of such an optimal s-curve by $\frac1L E(\ga,\gb)$, it is also shown
that $E(\ga,\gb)$ depends continuously on $(\ga,\gb)$. In the constructive proof of existence,
all optimal s-curves are described, but uniqueness is only proved in the case when the optimal curve is a c-curve, but not a u-turn.
An optimal s-curve is of {\bf Form 1} (resp. {\bf Form 2}) if it does not (resp. does) contain a u-turn. Optimal s-curves
of Form 1 are either line segments or segments of rectangular elastica (see Fig.\ \figsixteen\ (a)) while those of Form 2 
(see Fig.\ \figsixteen\ (b))
contain a u-turn of rectangular elastica along with, possibly, line segments and a c-curve of rectangular elastica.

Elastic splines were computed in a computer program {\it Curve Ensemble}, written in conjunction with \cite{\JJ}, and it was
observed that the fairness of elastic splines can be significantly degraded when pieces of Form 2 arise. As a remedy, it was suggested
that elastic splines be further restricted by requiring that chord angles of pieces satisfy
$$
\abs\ga,\abs\gb\leq\frac\pi2.
\tag\Elabel\defrestricted
$$ 
This additional restriction, which is stronger than (\scurveineq),
also greatly simplifies the numerical computation and theoretical development, and for these reasons, we have elected to 
adopt this restriction and so define $$\Cal A_{\pi/2}(P_1,P_2,\ldots,P_m)$$ to be the set of curves in
$\Cal A(P_1,P_2,\ldots,P_m)$ whose pieces have chord angles satisfying (\defrestricted). 
Curves in $\Cal A_{\pi/2}(P_1,P_2,\ldots,P_m)$
with minimal bending energy are called {\bf restricted elastic splines}.

In Section 3, we show that if (\defrestricted) holds and  $(\ga,\gb)\not\in\set{(\pi/2,-\pi/2),(-\pi/2,\pi/2)}$,
then the optimal s-curve from $P$ to $Q$, with chord angles $(\ga,\gb)$, is unique and of Form 1.
The omitted cases correspond to u-turns (see Fig.\ \figseventeen\ (a)) which fail to be unique only because one
can always extend a u-turn with line segments without affecting optimality. Nevertheless, the u-turn of rectangular
elastica (see Fig.\ \figseventeen\ (b)) is the unique $C^\infty$ optimal s-curve when  $(\ga,\gb)\in\set{(\pi/2,-\pi/2),(-\pi/2,\pi/2)}$.
We mention, belatedly, that the optimality of the u-turn of rectangular elastica was first proved by
Linn\'er and Jerome \cite{\LinnerJerome}.
\newline
\epsfysize=1truecm \epsffile{space.eps}\hskip1cm \epsfysize=3truecm \epsffile{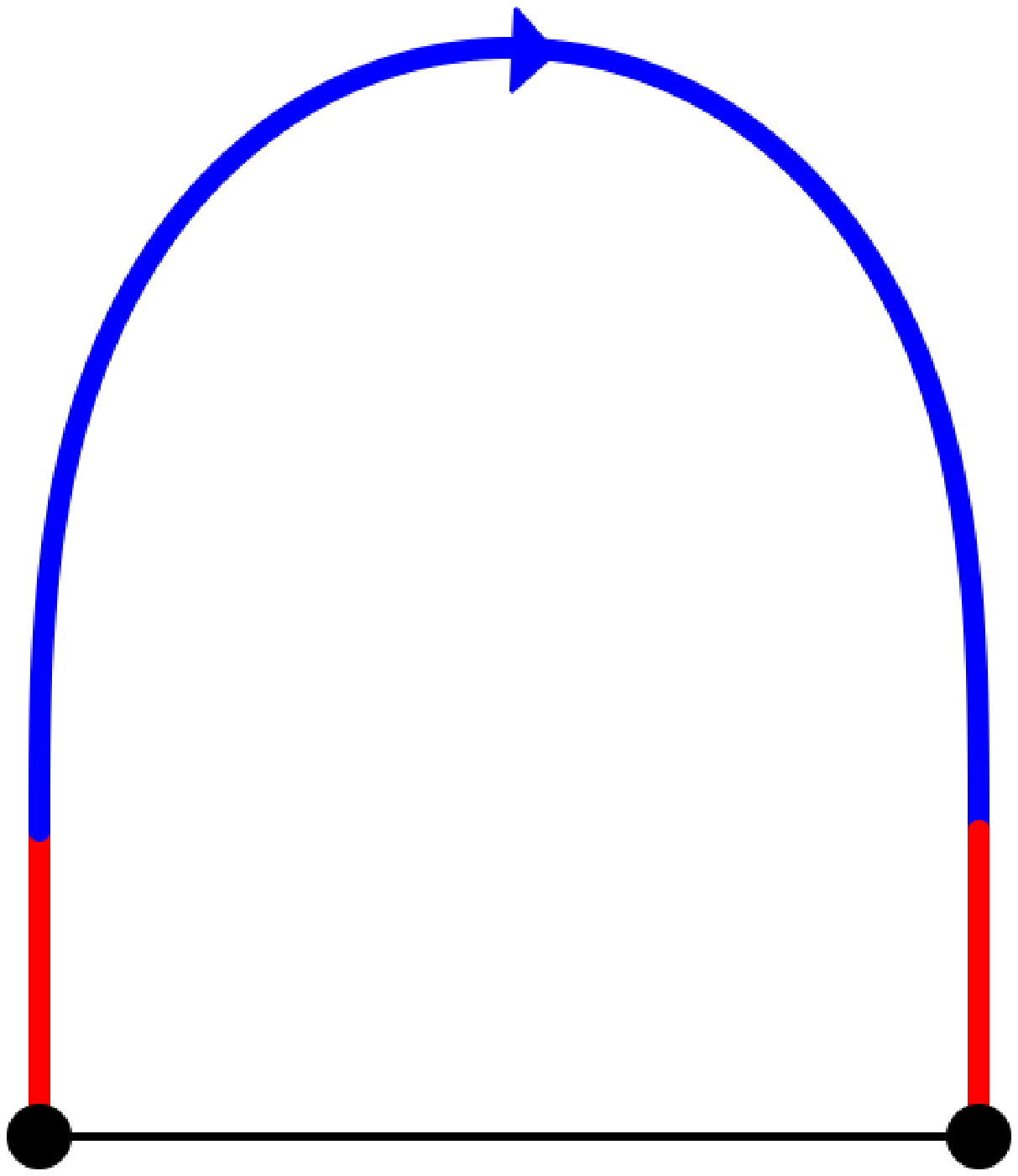}\hskip2.5cm
\epsfysize=3truecm \epsffile{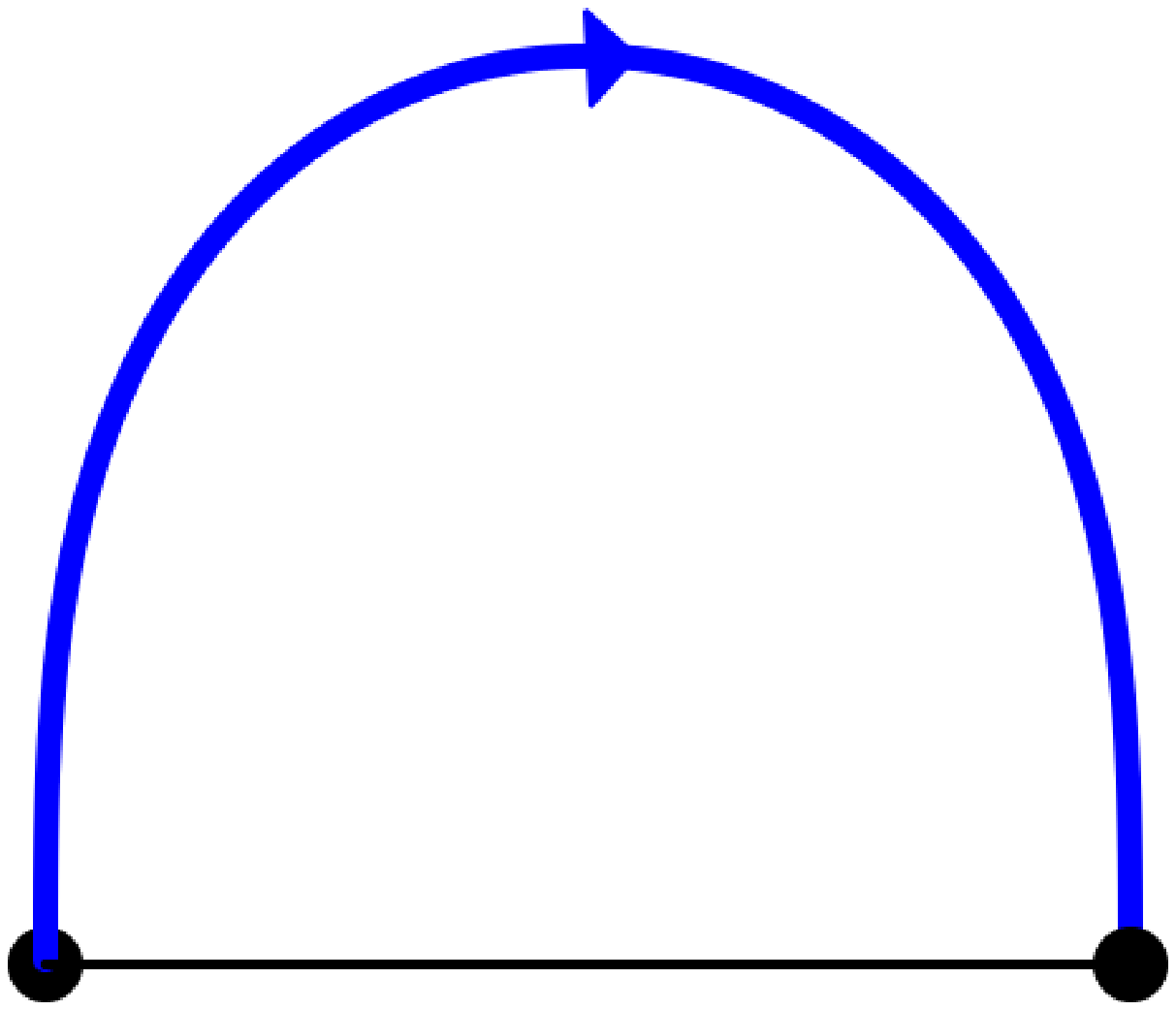}\newline
{\bf Fig.\ \Flabel\figseventeen} \ \ \ {\bf (a)} optimal u-turn \ \ \ \ \ \ \ \ \ \ \ \ \ \ \ \ 
{\bf (b)} u-turn of rectangular elastica.\newline
With unicity of optimal s-curves in hand, we can then appeal to the framework developed in \cite{\JJ} for assistance in
proving existence and $G^2$-regularity of restricted elastic splines. 
The following will be proved in Section \SectionRES.
\proclaim{Proposition \Tlabel\Fopt} The set $\Cal A_{\pi/2}(P_1,P_2,\ldots,P_m)$ contains a curve $F_{opt}$
with minimal bending energy. Moreover, if $F\in \Cal A_{\pi/2}(P_1,P_2,\ldots,P_m)$ has minimal bending energy,
then each piece of $F$ is $G^2$.
\endproclaim
\remark{Remark}
When discussing geometric curves, the notions of geometric regularity, $G^1$ and $G^2$, are preferred over the more
familiar notions of parametric regularity, $C^1$ and $C^2$. A curve $F$ has $G^1$ regularity if its unit tangent direction changes
continuously with respect to arclength and it has $G^2$ regularity if, additionally, its signed curvature changes
continuously with arclength. By our definition of curve (given at the outset), all curves are $G^1$, but not necessarily
$G^2$.
\endremark
The main concern of the present contribution is to identify conditions under which a restricted elastic spline
$F_{opt}$ will be globally $G^2$.
This direction of inquiry is motivated by a result
of Lee \& Forsyth \cite{\LeeForsyth} (see also Brunnett \cite{\Brunnett}) which says that if an interpolating
curve $F$ has bending energy which is locally minimal (i.e., minimal among all `nearby' interpolating curves),
then $F$ is globally $G^2$.
The proofs in \cite{\LeeForsyth}
and \cite{\Brunnett} employ variational calculus, but we prefer the constructive approach of \cite{\JJ} for
its clarity and generality.
We now explain our results on $G^2$ regularity assuming that  
$F_{opt}$ is a curve in $\Cal A_{\pi/2}(P_1,P_2,\ldots,P_m)$ having minimal bending energy. 
Note that it does not follow from Proposition \Fopt\ that $F_{opt}$ is globally $G^2$ because it is possible
for the signed curvature to have jump discontinuities across the interior nodes $P_2,P_3,\ldots,P_{m-1}$.
The following is a consequence of Theorem \refinedcondCtwo.
%
\proclaim{Corollary \Tlabel\condCtwo} If the chord angles at interior nodes are all (strictly) less than
$\frac\pi2$ in magnitude, then $F_{opt}$ is globally $G^2$.
\endproclaim
Proposition \Fopt\ and Corollary \condCtwo\ are analogous to results of Jerome and Fisher \cite{\JeromeOne, 
\JeromeTwo, \FisherJerome}
in that first additional constraints are imposed in order to ensure existence of an optimal curve, and then it is
shown that if these additional constraints are inactive, the optimal curve is globally
$G^2$ and its pieces are segments of rectangular elastica. These results are a good start, but they are
not entirely satisfying because they shed no light on whether one can expect the added constraints to be inactive.

Our experience using the program
{\it Curve Ensemble}  is that the hypothesis of Corollary \condCtwo\ holds
when the interpolation points $\set{P_j}$ impose only mild changes in direction. 
This vague idea can be quantified in
terms of the {\bf stencil angles} $\set{\psi_j}$ (see Fig.\ \figtwelve),  defined by 
$$
\psi_j:=\arg \frac{P_{j+1}-P_j}{P_j-P_{j-1}},\quad j=2,3,\ldots,m-1.
$$
\epsfysize=2.5truecm \epsffile{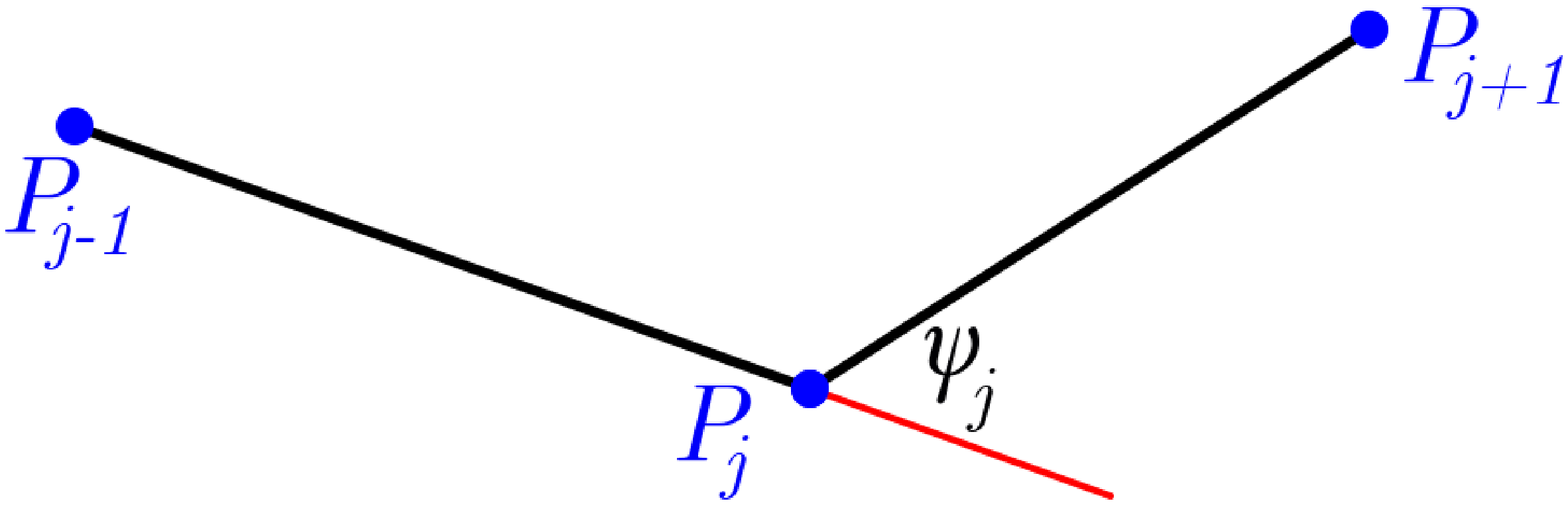}
\hskip1.5truecm \epsfysize=1.5truecm \epsffile{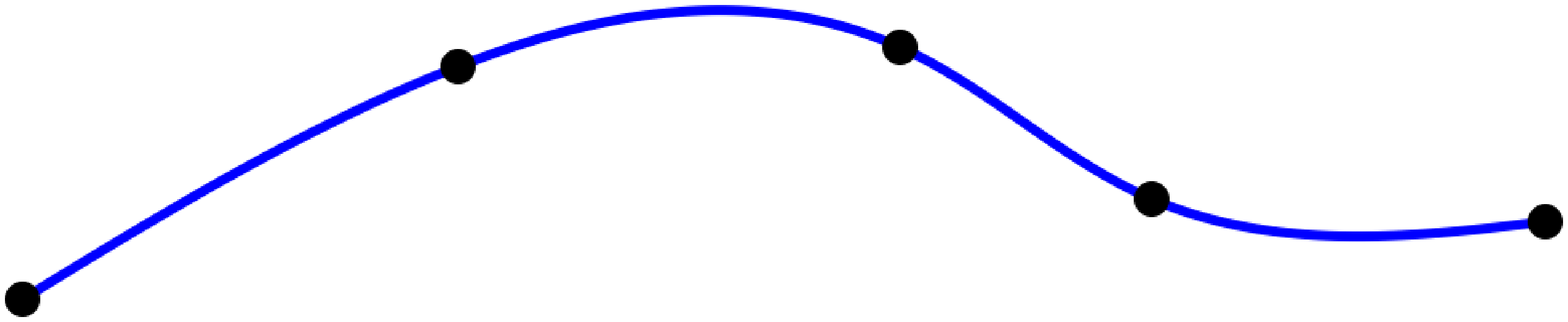} \newline
{\bf Fig.\ \Flabel\figtwelve} \ the stencil angle $\psi_j$
\hskip2truecm {\bf Fig.\ \Flabel\figthirteen} \ a globally $G^2$ restricted elastic spline\newline
The following is a consequence of Theorem \nothingyet.
\proclaim{Corollary \Tlabel\Ctwo} Let $\Psi$ ($\approx 37^\circ$) be the positive angle defined in
(\Psidefine). If the stencil angles satisfy
$\abs{\psi_j}< \Psi$ for $j=2,3,\ldots,m-1$, then the hypothesis of Corollary \condCtwo\ holds and
consequently $F_{opt}$ is globally $G^2$.
\endproclaim
For example, the stencil angles in Fig.\ \figthirteen\ are all less than $\Psi$ and therefore it follows from
Corollary \Ctwo\ that the shown restricted elastic spline is globally $G^2$.

%

An outline of the remainder of the paper is as follows. In Section 2, we summarize some notation from
\cite{\BJone} which is needed here, and then in Section 3, as mentioned above, we address the unicity
of optimal s-curves.
%
The proofs of our results on $G^2$ regularity are complicated by the fact that they are obtained by combining a variety of
related results, and so, for the sake of readability, we will `prove' these results
in Section \SectionRES, but leave the proofs of two key identities, namely (\consistent) and (\horseman),  
to later sections. Furthermore,
the proofs given in Section \SectionRES\ make essential use of the framework established in \cite{\JJ}, and so
Section \SectionRES\ begins by defining a {\it basic curve method}, called Restricted Elastic Splines,
which fits into the framework of \cite{\JJ}. Identities (\consistent) and (\horseman) are proved
in sections 7 and 8, but these proofs require a great deal of preparation (sections \SectionQ,\SectionQunique) relating
to the chord angles of parametrically defined segments of rectangular elastica. In addition to supporting the proofs
in sections 7,8, the preparations done in sections 5,6 are also useful
in the efficient numerical computation of restricted elastic splines.
%
\Section{Summary of Notation}
%
The present contribution uses the same notation as in \cite{\BJone}; we summarize it here. As mentioned above,
a curve is a function $f:[a,b]\to\BbC$ whose derivative $f'$ is absolutely continuous and non-vanishing.  The
{\bf bending energy} of $f$ is defined by
$$
\norm{f}^2 := \frac14\int_0^{\Cal L} \gk^2\,ds,
$$
where $\Cal L$ denotes the arclength of $f$ and $\gk$ its signed curvature (the unusual factor $1/4$ is used to simplify some
formulae related to rectangular elastica). Let $g:[c,d]\to\BbC$ be another curve.  We say that $f$ and $g$ are
{\bf equivalent} if they have the same arclength parametrizations.  They are {\bf directly similar} if there exists
a linear transformation $T(z)=c_1z+c_2$ ($c_1,c_2\in\BbC$) such that $f$ and $T\circ g$ are equivalent; if
$|c_1|=1$, they are called {\bf directly congruent}. The notions of {\bf similar} and {\bf congruent} are the same
except that $T$ is allowed to have the form $T(z)=c_1\overline z + c_2$, where $\overline z$ denotes the complex
conjugate of $z$.

As mentioned earlier,
we call $f$ an {\bf s-curve} if it first turns monotonically at most $180^\circ$ in one direction and then turns at most
$180^\circ$ in the opposite direction. 
An s-curve which turns in only one direction is called a {\bf c-curve}
and a c-curve which turns a full $180^\circ$ is called a {\bf u-turn}. 
A non-degenerate s-curve is called a {\bf left-right s-curve} if it first turns clockwise and then turns counter-clockwise;
otherwise it is called a {\bf right-left s-curve}.
S-curves are often associated with a
geometric Hermite interpolation problem, and so to facilitate this we employ the unit tangent vectors 
$u=(f(a),f'(a)/|f'(a)|)$ and $v=(f(b),f'(b)/|f'(b)|)$ to say that $f$ {\bf connects} $u$ to $v$.
If $g:[c,d]\to\BbC$ is a curve satisfying $(g(c),g'(c)/|g'(c)|)=(f(b),f'(b)/|f'(b)|)$, then $f\splice g$ denotes
the concatenated curve which, for the sake of clarity, is assumed to have the arclength parametrization.
Most of the s-curves which we will encounter are segments of rectangular elastica; our preferred
parametrization is $R(t)= \sin t + i\, \xi(t)$, where $\xi(t)$ is defined by
$\dsize \frac{d\x}{dt}=\frac{\sin^2 t}{\sqrt{1+\sin^2 t}}$, $\xi(0)=0$.
One easily verifies that $\x$ is odd and satisfies
$\xi(t+\pi)=d+\xi(t)$, where $d:=\xi(\pi)$. Since the sine function is odd and $2\pi$-periodic, we conclude that
$R(t)$ is odd and satisfies $R(t+2\pi)=i\, 2d+R(t)$.
For later reference, we mention the following.
$$\align
\abs{R'(t)}=\frac{1}{\sqrt{1+\sin^2 t}}&,\qquad
\frac{R'(t)}{\abs{R'(t)}}=\cos t \sqrt{1+\sin^2 t}+ i\, \sin^2 t ,\qquad \gk(t) = 2\sin t,\\
\norm{R_{[a,b]}}^2 &=\frac14\int_a^b \gk(t)^2|R'(t)|\,dt= \xi(b)-\xi(a),\\
\endalign
$$
where $R_{[a,b]}$ denotes the restriction of $R$ to the interval $[a,b]$.
%
\Section{Unicity of optimal s-curves}
%
Let $\ga,\gb\in(-\pi,\pi]$ and set $u=(0,e^{i\ga})$ and $v=(1,e^{i\gb})$. The set $S(\ga,\gb)$, defined
to be the set of all s-curves connecting $u$ to $v$, was intensely studied in \cite{\BJone}, and it is
easy to verify that $S(\ga,\gb)$ is non-empty if and only if
$(\ga,\gb)\in\Cal F$, where
$$
\Cal F:=\set{(\ga,\gb):\abs\ga,\abs\gb<\pi\text{ and }\abs{\ga-\gb}\leq \pi}.
$$
It is shown in \cite{\BJone} that if $S(\ga,\gb)$ is non-empty, then $S(\ga,\gb)$ contains a curve with
minimal bending energy; that is, there exists a curve $f_{opt}\in S(\ga,\gb)$ such that
$\norm{f_{opt}}^2\leq \norm{f}^2$ for all $f\in S(\ga,\gb)$. The bending energy of $f_{opt}$ is denoted
$$
E(\ga,\gb):=\norm{f_{opt}}^2,\quad (\ga,\gb)\in\Cal F.
\tag\Elabel\defE$$
Let $\Sopt(\ga,\gb)$ denote the set of all arclength parameterized curves in $S(\ga,\gb)$ whose bending energy is minimal.
In \cite{\BJone}, every curve in $\Sopt(\ga,\gb)$ is `described', but uniqueness
is only established in a few cases. In the present section, we obtain uniqueness results (Theorem \sunique)
for the case when $(\ga,\gb)$ belongs to the square $[-\frac\pi2,\frac\pi2]^2$ (note that $[-\frac\pi2,\frac\pi2]^2$
is the largest square of the form $[-\gO,\gO]^2$ which is contained in $\Cal F$).

%
\proclaim{Theorem \Tlabel\sunique} Assume $(\ga,\gb)\in [-\frac\pi2,\frac\pi2]^2$.
Then $\Sopt(\ga,\gb)$ contains a unique $C^\infty$ curve $c_1(\ga,\gb)$. Moreover, the following hold.\newline
(i) If $\abs{\ga-\gb}<\pi$, then $\Sopt(\ga,\gb)=\set{c_1(\ga,\gb)}$.\newline
(ii) If $\abs{\ga-\gb}=\pi$, then every curve in $\Sopt(\ga,\gb)$ is $C^2$. \newline
(iii) If $(\ga,\gb)\neq (0,0)$, then there exist $t_1<t_2<t_1+2\pi$ such that $c_1(\ga,\gb)$ is directly similar to $R_{[t_1,t_2]}$.
\endproclaim

Since the bending energy of a curve is invariant under translations, rotations, reflections and reversals (of orientation), 
when proving items (i) and (ii),
we can additionally assume, without loss of generality, that $\ga\geq\abs\gb$. This reduction is also valid for item (iii)
since $R_{[t_1+\pi,t_2+\pi]}$ is directly congruent to reflections of $R_{[t_1,t_2]}$ and $R_{[-t_2,-t_1]}$ is directly congruent 
to the reversal of $R_{[t_1,t_2]}$.
Our proof of Theorem \sunique\ uses some definitions and results from \cite{\BJone} which are posed assuming
$$
\ga\in(0,\pi),\quad \abs\gb\leq\ga,\quad \gb>\ga-\pi.
\tag\Elabel\canonical$$

In \cite{\BJone, section 5}, the following functions of $\gga\in\gG:=[\ga-\pi,\gb]\cap (-\infty,0)$ are introduced:
$$\align
y_1&:=y_1(\gga):=\frac12\int_0^{\ga-\gga}\sqrt{\sin\tau}\dd \tau \\
y_2&:=y_2(\gga):=\frac12\int_0^{\gb-\gga}\sqrt{\sin\tau}\dd \tau \\
G(\gga)&:=\frac1{-\sin\gga}(y_1+y_2)^2\\
\gs(\gga)&:=\cos\gga +\frac{\sin\gga}{y_1+y_2}(\sqrt{\sin(\ga-\gga)}+\sqrt{\sin(\gb-\gga)})\\
q(\gga)&:=\frac{-\sin\gga}{y_1+y_2}\\
\endalign
$$

If $\gb\geq 0$, then all curves in $S(\ga,\gb)$ are non-degenerate right-left s-curves and so the same is true of
$\Sopt(\ga,\gb)$.  In contrast, if $\gb<0$ then $S(\ga,\gb)$
contains right c-curves as well as non-degenerate s-curves (both right-left and left-right). Nevertheless,
it turns out that the curves in $\Sopt(\ga,\gb)$ are all of the same flavor.
The discriminating factor, when $\gb<0$, is the quantity $\gs(\gb)$:\newline
If $\gs(\gb)>0$, then all curves
in $\Sopt(\ga,\gb)$ are non-degenerate right-left s-curves, while if $\gs(\gb)\leq 0$, then the unique curve
in  $\Sopt(\ga,\gb)$ is a right c-curve.

Regarding the latter case, we have the following which follows from \cite{\BJone, Theorem 6.2}.
\proclaim{Theorem \Tlabel\sixtwo} Let (\canonical) be in force and assume that $\gb<0$ and $\gs(\gb)\leq 0$.
Then there exists a unique $C^\infty$ curve $c(\ga,\gb)$ such that  $\Sopt(\ga,\gb)=\set{c(\ga,\gb)}$. Furthermore,
there exist $-\pi<t_1<t_2\leq 0$ such that $c(\ga,\gb)$ is directly similar to $R_{[t_1,t_2]}$.
\endproclaim
%
The following lemma and proposition are consequences of
\cite{\BJone, Lemma 5.11} and \cite{\BJone, Corollary 5.12 and Remark 5.10}, respectively.
\proclaim{Lemma \Tlabel\randang} Assume (\canonical). The function $G:\gG\to (0,\infty)$ is continuously differentiable, has a minimum
value $G_{\text{min}}$, and satisfies $\dsize \frac d{d\gga} G(\gga) = \frac 1{q(\gga)^2}\gs(\gga)$ for all $\gga\in\gG$.
\endproclaim
\proclaim{Proposition \Tlabel\fivetwelve } Let (\canonical) be in force and in case $\gb<0$, assume $\gs(\gb)>0$.
Suppose that there exists $\wh\gga\in\gG$, with $\wh\gga>\ga-\pi$, such that $G$ is uniquely minimized at $\wh\gga$
(i.e., $G(\wh\gga)=G_{\text{min}}$ and
$G(\gga)>G_{\text{min}}$ for all $\gga\in\gG\bs\set{\wh\gga}$).  Then $\gs(\wh\gga)=0$ and
there exists a unique $C^\infty$ curve $c(\ga,\gb)$
such that  $\Sopt(\ga,\gb)=\set{c(\ga,\gb)}$. Moreover, $E(\ga,\gb)=G_{\text{min}}$ and
$c(\ga,\gb)$ is directly similar to $R_{[t_1,t_2]}$,
where $-\pi<t_1<0<t_2<\pi$ are uniquely determined by $\arg R'(t_1)=\ga-\wh\gga$ and $\arg R'(t_2)=\gb-\wh\gga$.
\endproclaim
\remark{Remark} That the above conditions $\arg R'(t_1)=\ga-\wh\gga$ and $\arg R'(t_2)=\gb-\wh\gga$
do determine  $-\pi<t_1<0<t_2<\pi$ uniquely can be verified by first noting that $\arg R'(t)$ decreases
continuously from $\pi$ to $0$, as $t$ runs from $-\pi$ to $0$, and then increases continuously back up to $\pi$
as $t$ runs from $0$ to $\pi$. Now, since (\canonical) holds, $\wh\gga\in\gG$, and $\wh\gga>\ga-\pi$, it follows that
$0<\ga-\wh\gga<\pi$ and $0\leq\gb-\wh\gga<\pi$. What remains is to show that $0<\gb-\wh\gga$. If $\gb\geq 0$,
then  $0<\gb-\wh\gga$ is clear since $\wh\gga<0$. If $\gb<0$, then we cannot have $\gb=\wh\gga$ because
$\gs(\gb)>0$ while $\gs(\wh\gga)=0$; therefore,  $0<\gb-\wh\gga$.
\endremark
We now begin the proof of Theorem \sunique, and, as mentioned above, it suffices to prove the theorem in the
canonical case when $(\ga,\gb)\in[-\pi/2,\pi/2]^2$ satisfy $\ga\geq\abs\gb$. We begin with two specific cases.
\demo{Proof of Theorem \sunique\ in case $(\ga,\gb)=(0,0)$}
In this case, it is easy to verify that $\Sopt(0,0)$ contains only the line segment from $0$ to $1$.
\qed\enddemo
\demo{Proof of Theorem \sunique\ in case  $(\ga,\gb)=(\frac\pi2,-\frac\pi2)$}
By definition of an s-curve, every curve in $S(\frac\pi2,-\frac\pi2)$ is a right u-turn
and it is shown in \cite{\BJone, sections 3,4} that every curve in  $\Sopt(\frac\pi2,-\frac\pi2)$ is either directly similar
to $R_{[-\pi,0]}$ or else equals $[0,iq]\splice f \splice [1+iq,1]$ where $q>0$ and $f$ is directly similar to  $R_{[-\pi,0]}$
(here, $[0,iq]$ and $[1+iq,1]$ denote line segments).
Among these, the only curve which is $C^\infty$ is the first one, and therefore
$\Sopt(\frac\pi2,-\frac\pi2)\cap C^\infty = \set{c(\frac\pi2,-\frac\pi2)}$,
where $c(\frac\pi2,-\frac\pi2)$ is the arclength parameterized curve in $S(\frac\pi2,-\frac\pi2)$ which is directly similar to $R_{[-\pi,0]}$. Since
the signed curvature of  $R_{[-\pi,0]}$ vanishes at the endpoints, it follows that all curves in $\Sopt(\frac\pi2,-\frac\pi2)$ are $C^2$,
which proves item (ii).
\qed\enddemo
Having proved Theorem \sunique\ in these two specific cases, we proceed assuming that
$$
\ga\in(0,\pi/2],\quad \abs{\gb}\leq\ga,\quad \gb>-\pi/2,
\tag\Elabel\rescanonical $$
and we note that (\rescanonical) implies (\canonical).
\proclaim{Lemma \Tlabel\wasabi} Assume (\rescanonical) and let $\gga\in\gG$. The following hold.\newline
(i) If $(\ga,\gb)\neq (\pi/2,\pi/2)$ and $\gga\leq -\pi/2$ , then $\gs(\gga)<0$.\newline
(ii) If $(\ga,\gb) = (\pi/2,\pi/2)$, then $\gs(-\pi/2)=0$ but there exists $\gep>0$ such that
$\gs(\gga)<0$ for all $\gga\in(-\pi/2,-\pi/2+\gep]$.\newline
(iii) If  $\gs(\gga)=0$ and $-\pi/2<\gga<\gb$, then $\gs'(\gga)>0$.
\endproclaim
\demo{Proof} We first note that $y_1>0$, $y_2\geq 0$, and $\gs$ is continuous on $\gG$.
And since both $\sin(\ga-\gga)$ and $\sin(\gb-\gga)$ are positive for $\gga$ in the interior $\gG^{o}:=(\ga-\pi,\gb)\cap(-\infty,0)$,
it follows that $\gs$ is $C^1$ on $\gG^{o}$.
Defining $H(\gga):=y_1+y_2$, $\gga\in\gG$, we have
$$
H'(\gga)=-\frac12\pr{\sqrt{\sin(\ga-\gga)}+\sqrt{\sin(\gb-\gga)}},\quad
H''(\gga)=\frac14\pr{\frac{\cos(\ga-\gga)}{\sqrt{\sin(\ga-\gga)}}
+\frac{\cos(\gb-\gga)}{\sqrt{\sin(\gb-\gga)}}}
$$
and it follows that $H$ is $C^1$ on $\gG$ and $C^2$ on $\gG^{o}$. Moreover, we note that $H$ is positive on $\gG$, while
$H'(\gga)<0$ for all $\gga\in\gG$ except when $-\gga=\ga=\gb=\pi/2$. We can express $\gs$ in terms of $H$ as
$$
\gs(\gga) H(\gga) = \cos\gga H(\gga) - 2\sin\gga H'(\gga),\quad \gga\in\gG,
\tag\Elabel\dayeight$$
and then differentiation yields
$$
\gs'(\gga)H(\gga)+\gs(\gga)H'(\gga) = -\sin\gga H(\gga) -\cos\gga H'(\gga) - 2\sin\gga H''(\gga),\quad \gga\in\gG^{o}.
\tag\Elabel\daynine$$
We first prove (i): Assume $(\ga,\gb)\neq (\pi/2,\pi/2)$ and $\gga\leq -\pi/2$. Then $H(\gga)>0$, $\cos\gga\leq 0$,
$H'(\gga)<0$ and $\sin\gga<0$, and it follows easily from (\dayeight) that $\gs(\gga)<0$.\newline
We next prove (ii): Assume $(\ga,\gb) = (\pi/2,\pi/2)$. Then $\gG=[-\pi/2,0)$ and it is clear from the definition of
$\gs$ that $\gs(-\pi/2)=0$. In order to prove (ii), it suffices to show that
$\gs'(\gga)\to -\infty$ as $\gga\to -\pi/2^+$. Now, $H(-\pi/2)>0$, $H'(-\pi/2)=0$, but $H''(\gga)\to -\infty$ as
$\gga\to -\pi/2^+$. It therefore follows from (\daynine) that $\gs'(\gga)\to -\infty$ as $\gga\to -\pi/2^+$.\newline
Lastly, we prove (iii): Assume  $\gs(\gga)=0$ and $-\pi/2<\gga<\gb$. Since $\gs(\gga)=0$ and $\cos\gga>0$, it follows
from the definition of $\gs$ that $y_1+y_2 = -\tan\gga\pr{\sqrt{\sin(\ga-\gga)}+\sqrt{\sin(\gb-\gga)}}$; that is,
$H(\gga)=2\tan\gga H'(\gga)$. Substituting this into (\daynine) then yields
$$\align
\gs'(\gga) H(\gga) &= - \sin\gga \pr{2\tan\gga H'(\gga)+2H''(\gga)} - \cos\gga H'(\gga) \\
                   &= - \frac12\tan\gga \pr{4\sin\gga H'(\gga)+4\cos\gga H''(\gga)} - \cos\gga H'(\gga) \endalign
$$
Since $H(\gga)$, $-H'(\gga)$, $-\sin\gga$, and $ \cos\gga$ are positive, in order to prove that $\gs'(\gga)>0$, it suffices to
show that $4\sin\gga H'(\gga)+4\cos\gga H''(\gga)$ is nonnegative. Using the above formulations for $H'(\gga)$ and
$H''(\gga)$ and the identity $\cos(x+y)=\cos x \cos y - \sin x \sin y$, it is easy to verify that
$2\sin\gga H'(\gga)+4\cos\gga H''(\gga)=\frac{\cos\ga}{\sqrt{\sin(\ga-\gga)}} + \frac{\cos\gb}{\sqrt{\sin(\gb-\gga)}}$.
Hence, $4\sin\gga H'(\gga)+4\cos\gga H''(\gga)= 2\sin\gga H'(\gga) +
\frac{\cos\ga}{\sqrt{\sin(\ga-\gga)}} + \frac{\cos\gb}{\sqrt{\sin(\gb-\gga)}} > 0$.
\qed\enddemo
\demo{Proof of Theorem \sunique\ in case (\rescanonical) holds and $\gb<0$}
If $\gs(\gb)\leq 0$, then Theorem \sunique\ is an immediate consequence of Theorem \sixtwo; so assume that $\gs(\gb)>0$.
Note that $\gG=[\ga-\pi,\gb]$ and, by Lemma \wasabi\ (i), $\gs(\gga)<0$ for all $\gga\in[\ga-\pi,-\pi/2]$. Since $\gs$
is continuous and $\gs(\gb)>0$, it follows that there exists $\wh\gga\in(-\pi/2,\gb)$ such that $\gs(\wh\gga)=0$. It
follows from Lemma \wasabi\ (iii) that $\wh\gga$ is the only $\gga\in(-\pi/2,\gb)$ where $\gs$ vanishes. Therefore,
$\gs(\gga)<0$ for $\gga\in[\ga-\pi,\wh\gga)$ and $\gs(\gga)>0$ for $\gga\in(\wh\gga,\gb]$. It now follows from Lemma
\randang\ that $G$ is uniquely minimized at $\wh\gga$ and we obtain Theorem \sunique\ as a consequence of Proposition \fivetwelve.
\qed\enddemo
\demo{Proof of Theorem \sunique\ in case (\rescanonical) holds and $\gb\geq 0$}
Note that $\gG=[\ga-\pi,0)$. It follows from Lemma \wasabi\ (i) and (ii), and the continuity of $\gs$, that there exists
$\gep>0$ such that $\gs(\gga)<0$ for all $\gga\in(\ga-\pi,-\pi/2+\gep]$. From the definition of $\gs$, it is clear that
$\lim_{\gga\to 0^-}\gs(\gga)=1$, and hence there exists $\wh\gga\in(-\pi/2+\gep,0)$ such that $\gs(\wh\gga)=0$. As in the previous case,
it follows from Lemma \randang\ that $G$ is uniquely minimized at $\wh\gga$ and we obtain Theorem \sunique\ as a consequence
of Proposition \fivetwelve.
\qed\enddemo
This completes the proof of Theorem \sunique.
%
\Section{The Restricted Elastic Spline and Proofs of Main Results}
%
Although written specifically for s-curves which connect $u=(0,e^{i\ga})$ to $v=(1,e^{i\gb})$, Theorem \sunique\ easily
extends to general configurations $(u,v)$.  To see this, let $u=(P_1,d_1)$ and $v=(P_2,d_2)$ be two unit tangent vectors 
with distinct base points $P_1\neq P_2$.  The chord angles $(\ga,\gb)$ determined by $(u,v)$ are
$\ga=\arg\frac{d_1}{P_2-P_1}$ and  $\gb=\arg\frac{d_2}{P_2-P_1}$.  With $S(u,v)$ denoting the set of s-curves which connect $u$ to $v$,
and defining $T(z):=(P_2-P_1)z+P_1$, we see that
$S(u,v)$ is in one-to-one correspondence with $S(\ga,\gb)$ (defined in Section 3):
$f\in S(\ga,\gb)$ if and only if $T\circ f \in S(u,v)$. Moreover, with $L:=|P_2-P_1|$, we have
$\norm{f}^2 = \frac1L \norm{T\circ f}^2$. Now, let us assume that $|\ga|,|\gb|\leq \pi/2$ and let $c_1(\ga,\gb)$ be the optimal 
arc-length parametrized curve described in
Theorem \sunique. Then $T\circ c_1(\ga,\gb)$ is an optimal curve in $S(u,v)$ having constant speed $L$ (not necessarily $1$), 
and so we define
$c(u,v)$ to be the arclength parametrized curve which is equivalent to $T\circ c_1(\ga,\gb)$. 
With $\Sopt(u,v)$ denoting the set of arclength parametrized curves in $S(u,v)$ having minimal bending energy,
Theorem \sunique\ translates immediately into the following.
\proclaim{Corollary \Tlabel\corsunique} Let $(u,v)$ be a configuration with chord angles $(\ga,\gb)\in [-\frac\pi2,\frac\pi2]^2$.
Then $c(u,v)$ is the unique $C^\infty$ curve in $\Sopt(u,v)$. Moreover, the following hold.\newline
(i) If $\abs{\ga-\gb}<\pi$, then $\Sopt(u,v)=\set{c(u,v)}$.\newline
(ii) If $\abs{\ga-\gb}=\pi$, then every curve in $\Sopt(u,v)$ is $C^2$. \newline
(iii)  $c(u,v)$ is directly similar to $c_1(\ga,\gb)$ and $\norm{c(u,v)}^2 = \frac1L \norm{c_1(\ga,\gb)}$.
\endproclaim

In the framework of \cite{\JJ}, the curves $\set{c(u,v)}$ are called {\bf basic curves} and the mapping 
$(u,v)\mapsto c(u,v)$ is called a {\bf basic curve method}. We define the {\bf energy} of basic curves to be
the bending energy. In \cite{\JJ}, it is assumed that the basic curve method and energy are translation and rotation
invariant, and this allows one's attention to be focused on
the (canonical) case where $u=(0,e^{i\ga})$ and $v=(L,e^{i\gb})$, $L>0$. The resulting basic curve and energy functional
are denoted
$c_L(\ga,\gb)$ and $E_L(\ga,\gb)$.  
In our setup, we have the two additional properties that
$c_L(\ga,\gb)$ is equivalent to $Lc_1(\ga,\gb)$ and $E_L(\ga,\gb) = \frac1L E_1(\ga,\gb)$, where the latter holds because
$$
E_L(\ga,\gb):=\norm{c_L(\ga,\gb)}^2 = \norm{Lc_1(\ga,\gb)}^2 =\frac1L\norm{c_1(\ga,\gb)}^2 = \frac1L E_1(\ga,\gb).
$$
In the language of \cite{\JJ}, we would say that the basic curve method is {\it scale invariant} and the energy functional
is {\it inversely proportional to scale}. This special case is addressed in detail in \cite{\JJ, sec. 3}, and it
allows us to focus our attention on the case $L=1$ where we have, for $(\ga,\gb)\in[-\frac\pi2,\frac\pi2]^2$, the optimal
curve $c_1(\ga,\gb)$ as described in Theorem \sunique\ and its energy $E_1(\ga,\gb)=\norm{c_1(\ga,\gb)}^2$.
Note that $E_1(\ga,\gb)=E(\ga,\gb)$ for  
$(\ga,\gb)\in [-\frac\pi2,\frac\pi2]^2$, where $E(\ga,\gb)$ is defined in (\defE).
The distinction between $E_1$ and $E$ is that the domain 
of $E_1$ is $[-\frac\pi2,\frac\pi2]^2$, while the domain of $E$ is the larger set $\Cal F$ (defined just above (\defE)).
In \cite{\BJone, sec. 7}, it is shown that $E$ is continuous on $\Cal F$ and it therefore follows that $E_1$ is
continuous on $[-\frac\pi2,\frac\pi2]^2$.  

The framework of \cite{\JJ} is concerned with the set $\wh{\Cal A}_{\pi/2}(P_1,P_2,\ldots,P_m)$ consisting of all
interpolating curves whose pieces are basic curves, and the energy of such an interpolating curve 
$\wh F=c(u_1,u_2)\splice c(u_2,u_3)\splice \cdots \splice c(u_{m-1},u_m)$ is define
to be the sum of the energies of its constituent basic curves:
Energy$(\wh F):=\sum_{j=1}^{m-1} \norm{c(u_j,u_{j+1})}^2=\norm{\wh F}^2$.
Note that $\wh{\Cal A}_{\pi/2}(P_1,P_2,\ldots,P_m)$ is a subset of $\Cal A_{\pi/2}(P_1,P_2,\ldots,P_m)$ and energy
in both sets is defined to be bending energy. Since $E_1$ is continuous on  $[-\pi/2,\pi/2]^2$,
it follows from \cite{\JJ, Th. 2.3} that there exists a curve in $\wh{\Cal A}_{\pi/2}(P_1,P_2,\ldots,P_m)$ with minimal bending
energy.
\remark{Remark} Whereas curves in $\Cal A(P_1,P_2,\ldots,P_m)$ with minimal bending energy are called 
{\bf elastic splines}, 
such curves in $\wh{\Cal A}_{\pi/2}(P_1,P_2,\ldots,P_m)$ are called {\bf restricted elastic splines}.
\endremark
The following lemma will be needed in our proof of Proposition \Fopt.
\proclaim{Lemma \Tlabel\projection} Given $F\in \Cal A_{\pi/2}(P_1,P_2,\ldots,P_m)$, let $u_1,u_2,\ldots,u_m$ be the
unit tangent vectors, with base-points $P_1,P_2,\ldots,P_m$, determined by $F$, and define
$$
\wh F:=c(u_1,u_2)\splice c(u_2,u_3)\splice \cdots \splice c(u_{m-1},u_m)\in  \wh{\Cal A}_{\pi/2}(P_1,P_2,\ldots,P_m).
$$
Then $\norm{F}^2\geq \norm{\wh F}^2$.
\endproclaim
The proof of the lemma is simply that the $j$-th piece of $F$ has bending energy at least $\norm{c(u_j,u_{j+1})}^2$
because it belongs to $S(u_j,u_{j+1})$ while  $c(u_j,u_{j+1})$ belongs to  $\Sopt(u_j,u_{j+1})$.
\demo{Proof of Proposition \Fopt} Since $\wh{\Cal A}_{\pi/2}(P_1,P_2,\ldots,P_m)$ is a subset of
$\Cal A_{\pi/2}(P_1,P_2,\ldots,P_m)$ and the former contains a curve with minimal bending energy, it follows
immediately from Lemma \projection\ that the latter contains a curve with minimal bending energy.  Now, assume
$F\in \Cal A_{\pi/2}(P_1,P_2,\ldots,P_m)$ has minimal bending energy, and let $\wh F$ be as in Lemma \projection.
Then $\norm{F}^2=\norm{\wh F}^2$ and we must have $\norm{F_{[t_j,t_{j+1}]}}^2=\norm{c(u_j,u_{j+1})}^2$ for $j=1,2,\ldots,m-1$.
Hence $F_{[t_j,t_{j+1}]}$ is equivalent to a curve in $\Sopt(u_j,u_{j+1})$ and it follows from Theorem \sunique\ that $F_{[t_j,t_{j+1}]}$ is $G^2$.
\qed\enddemo
%

The following definition is taken from \cite{\JJ, sec. 3}.
\definition{Definition} Let $F\in \wh{\Cal A}_{\pi/2}(P_1,P_2,\ldots,P_m)$ have minimal bending energy and let 
$(\ga_j,\gb_{j+1})$ be the chord angles of the the $j$-th piece of $F$.  We say that $F$ is
{\bf conditionally $G^2$} if $F$ is $G^2$ across $P_j$ whenever the two chord angles associated with node $P_j$ satisfy
$|\gb_j|,|\ga_j|<\pi/2$.
\enddefinition
Let $\gk_a(f)$ and $\gk_b(f)$ denote, respectively,
the initial and terminal signed curvatures of a curve $f$. 
The following result is an amalgam of \cite{\JJ, Th. 3.3 and Th. 3.5}.
%
\proclaim{Theorem \Tlabel\JJcond} If there exists a constant $\mu\in\Bbb R$ such that
$$
-\gk_a(c_1(\ga,\gb)) = \mu \frac\partial{\partial\ga} E_1(\ga,\gb)\quad\text{ and }\quad
\gk_b(c_1(\ga,\gb)) = \mu \frac\partial{\partial\gb} E_1(\ga,\gb)
\tag\Elabel\consistent
$$
for  all $(\ga,\gb)\in[-\frac\pi2,\frac\pi2]^2$ with $\abs{\ga-\gb}<\pi$,
then minimal energy curves in $\wh {\Cal A}_{\pi/2}(P_1,P_2,\ldots,P_m)$ are conditionally $G^2$.
\endproclaim
\remark{Remark} Although \cite{\JJ, Th. 3.3} is stated assuming that (\consistent) holds for all $(\ga,\gb)\in[-\frac\pi2,\frac\pi2]^2$,
the given proof remains valid under the weaker assumption that (\consistent) holds for all $(\ga,\gb)\in[-\frac\pi2,\frac\pi2]^2$
with $\abs{\ga-\gb}<\pi$.
\endremark
In the following sections, culminating in Theorem \nablaEQ, we will show that condition (\consistent) holds
with $\mu=2$. Together, Theorem \JJcond\ and Theorem \nablaEQ\ imply that minimal energy
curves in $\wh{\Cal A}_{\pi/2}(P_1,P_2,\ldots,P_m)$ are conditionally $G^2$; we can now prove that this also holds
for the larger set $\Cal A_{\pi/2}(P_1,P_2,\ldots,P_m)$.
\proclaim{Theorem \Tlabel\refinedcondCtwo} Let $F\in \Cal A_{\pi/2}(P_1,P_2,\ldots,P_m)$ have minimal bending energy.
Then $F$ is $G^2$ across $P_j$ (i.e., $\gk_b(F_{[t_{j-1},t_j]})=\gk_a(F_{[t_j,t_{j+1}]})$) 
whenever the two chord angles associated with node $P_j$ satisfy
$|\gb_j|,|\ga_j|<\pi/2$.
\endproclaim
\demo{Proof} Let $\wh F\in \wh{\Cal A}_{\pi/2}(P_1,P_2,\ldots,P_m)$ be as in Lemma \projection, and let
 $j\in\set{2,3,\ldots,m-1}$ be such that $|\gb_j|,|\ga_j|<\pi/2$. By Theorem \JJcond\ and Theorem \nablaEQ, $\wh F$
is $G^2$ across $P_j$.  The chord angles of the $j$-th piece of $F$ are $(\ga_j,\gb_{j+1})$ and since $\abs{\ga_j}<\pi/2$,
we must have $\abs{\gb_{j+1}-\ga_j}<\pi$ and it follows from Corollary \corsunique\ (i) that the $j$-th piece of
$F$ is equivalent to the $j$-th piece of $\wh F$. Similarly, since $\abs{\gb_j}<\pi/2$, the $(j-1)$-th piece
of $F$ is equivalent to the $(j-1)$-th piece of $\wh F$. We therefore have
$$
\gk_b(F_{[t_{j-1},t_j]})=\gk_b(\wh F_{[t_{j-1},t_j]})=\gk_a(\wh F_{[t_j,t_{j+1}]})=\gk_a(F_{[t_j,t_{j+1}]}).
$$
\qed\enddemo


For $t\in(0,\pi]$, let the chord angles of $R_{[0,t]}$ be denoted $\ga(0,t)$ and $\gb(0,t)$ (these definitions will
be extended in Section \SectionQ). 
In Corollary \cortbar, we show that there exists a unique $\tbar\in(0,\pi)$ such that $\gb(0,\tbar)=\frac\pi2$. 
Let $\Psi$ (see Fig.\ \figfifteen) denote the positive angle defined by
$$
\Psi:=\frac\pi2-\abs{\ga(0,\tbar)}. \qquad\qquad\qquad \text{\bf Fig.\ \Flabel\figfifteen}\epsfysize=4truecm \epsffile{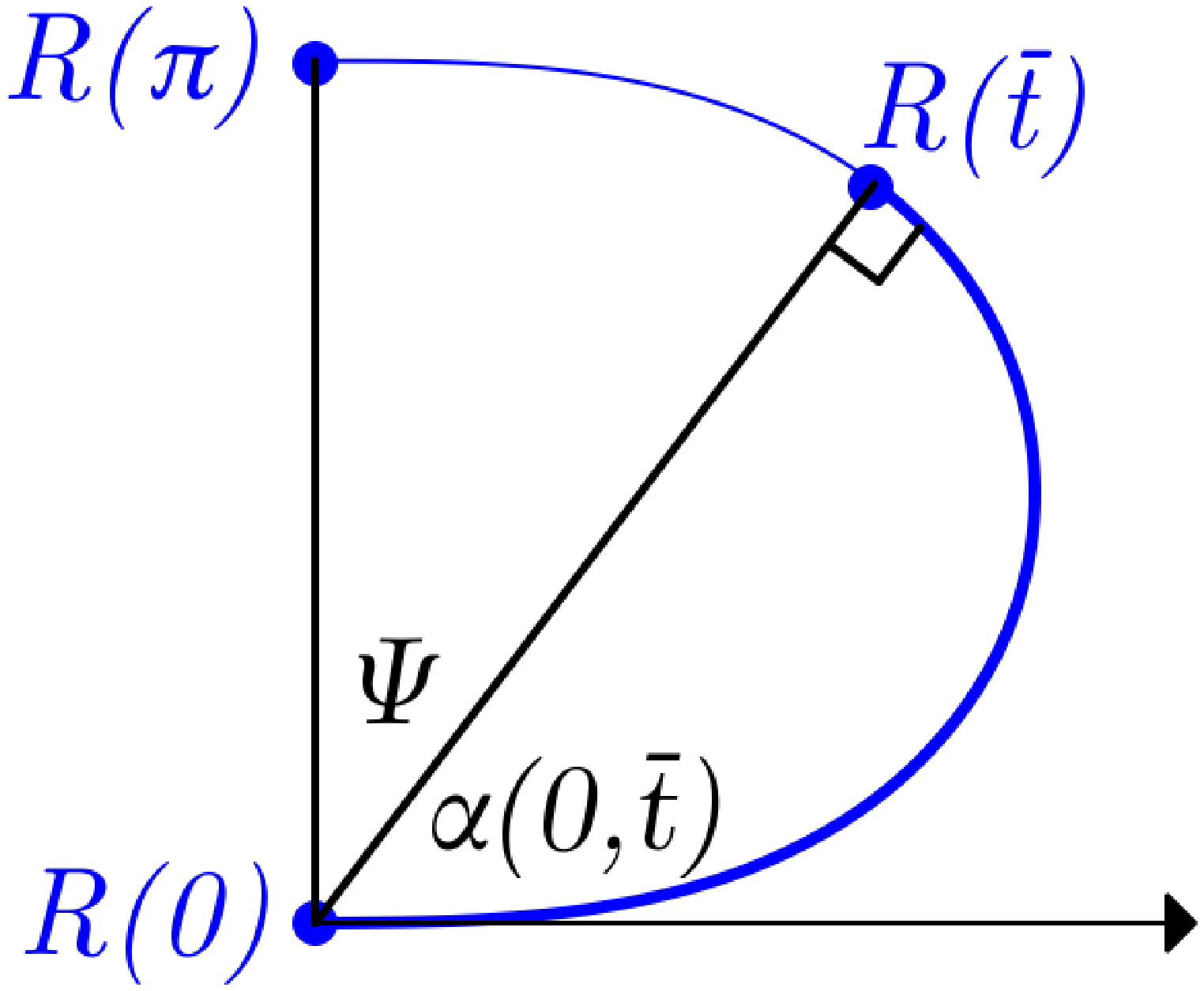}
\tag\Elabel\Psidefine$$
Our main result on $G^2$ regularity is obtained as a consequence of the following theorem which is
essentially \cite{\JJ, Theorem 5.1} but specialized to the present context.
\proclaim{Theorem \Tlabel\horsetheorem} Suppose that for every $\ga\in[-\frac\pi2,\frac\pi2]$ there exists $\gb_\ga^*$,
with $\abs{\gb_\ga^*}\leq \frac\pi2-\Psi$, such that
$$
\text{sign}\pr{\frac{\partial}{\partial\gb}E_1(\ga,\gb)}=\text{sign}(\gb - \gb_\ga^*)\text{ for all }
\gb\text{ satisfying }\abs\gb\leq\frac\pi2\text{ and }\abs{\gb-\ga}<\pi.\tag\Elabel\horseman$$
Let $\wh F\in\wh{\Cal A}_{\pi/2}(P_1,P_2,\ldots,P_m)$ be a curve with minimal bending energy. If $P_j$ 
is a point where the stencil angle $\psi_j$ satisfies $\abs{\psi_j}<\Psi$, then
the two chord angles associated with node $P_j$ satisfy
$|\gb_j|,|\ga_j|<\pi/2$ and, consequently, $\wh F$ is $G^2$ across node $P_j$.
\endproclaim
\demo{Proof} Employing the symmetry  $E_1(\ga,\gb)=E_1(\gb,\ga)$, conditions (i) and (ii) in the hypothesis of 
\cite{\JJ, Theorem 5.1} reduce simply to the single condition
$$
\text{sign}\pr{\frac{\partial}{\partial\gb}E_1(\ga,\gb)}=\text{sign}(\gb - \gb_\ga^*)\text{ for all }
\abs\gb\leq\frac\pi2.
\tag\Elabel\horsemanalt
$$
and therefore Theorem \horsetheorem, with (\horseman) replaced by (\horsemanalt), is an immediate consequence of 
\cite{\JJ, Theorem 5.1}. Note that the only distinction between (\horseman) and (\horsemanalt) is that (\horseman)
is mute when $(\ga,\gb)$ equals $(\pi/2,-\pi/2)$ or $(-\pi/2,\pi/2)$. 
With a slight modification (specifically: rather than showing that $f'(\gO)>0$ and 
$f'(\psi_2-\gO)<0$, one instead shows that there exists $\gep>0$ such that $f'(\gb)>0$ for $\gO-\gep<\gb<\gO$ and
$f'(\gb)<0$ for $\psi_2-\gO<\gb<\psi_2-\gO+\gep$), the proof of \cite{\JJ, Theorem 5.1} also proves
Theorem \horsetheorem. 
\qed\enddemo
\remark{Remark} The appearance of (\horseman), rather than (\horsemanalt), in Theorem \horsetheorem\ is simply
a consequence of the fact that $\frac{\partial}{\partial\gb}E_1(\ga,\gb)=0$ when $(\ga,\gb)$ equals $(\pi/2,-\pi/2)$ 
or $(-\pi/2,\pi/2)$. This distinction is not without consequence. In \cite{\JJ, Theorem 5.1}, the conclusion is obtained
when $\psi_i\leq\Psi$, while in Theorem \horsetheorem\ we require $\psi_i<\Psi$.
\endremark

In Section 8, we prove that (\horseman) holds and we therefore obtain the conclusion of Theorem \horsetheorem\ 
regarding minimal energy curves in $\wh{\Cal A}_{\pi/2}(P_1,P_2,\ldots,P_m)$. We will now show that the same holds
for the larger set ${\Cal A}_{\pi/2}(P_1,P_2,\ldots,P_m)$. 
\proclaim{Theorem \Tlabel\nothingyet} Let $F\in \Cal A_{\pi/2}(P_1,P_2,\ldots,P_m)$ have minimal bending energy.
Then $F$ is $G^2$ across $P_j$ whenever the stencil angle satisfies $\abs{\psi_j}<\Psi$.
\endproclaim
\demo{Proof} Let $\wh F\in \wh{\Cal A}_{\pi/2}(P_1,P_2,\ldots,P_m)$ be as in Lemma \projection, and let
$j\in\set{2,3,\ldots,m-1}$ be such that
$\abs{\psi_j}<\Psi$. It follows from Theorem \horsetheorem\ and Section 9 that the two chord angles at
node $P_j$ satisfy $\abs{\gb_j},\abs{\ga_j}<\frac\pi2$, and therefore, by Theorem \refinedcondCtwo, $F$
is $G^2$ across $P_j$.
\qed\enddemo
%
\Section {The chord angles of  $R_{[t_1,t_2]}$}
%
In this section and the next, we establish relations between the {\bf parameters} $(t_1,t_2)$, with $t_1<t_2$,
and the chord angles $(\ga,\gb)$ of the segment $R_{[t_1,t_2]}$ of rectangular elastica (defined in Section 2). 
Our primary purpose in this section is to prove Theorem \detq\ and Corollary \cortstar.

Recall from Section 2 that
the chord angles are given by 
$\ga:=\ga(t_1,t_2)=\arg \frac{R'(t_1)}{R(t_2)-R(t_1)}$ and $\gb:=\gb(t_1,t_2)=\arg \frac{R'(t_2)}{R(t_2)-R(t_1)}$. 
We mention that since $\xi(t)$ is   increasing, it follows that the chord angles $\ga(t_1,t_2)$ and
$\gb(t_1,t_2)$ never equal $\pi$ (i.e., the branch cut in the definition of $\arg$ is never crossed).
\newline
{\bf Fig.\ \Flabel\figseven} \ notation for  $R_{[t_1,t_2]}$
\hskip1truecm \epsfysize=5truecm \epsffile{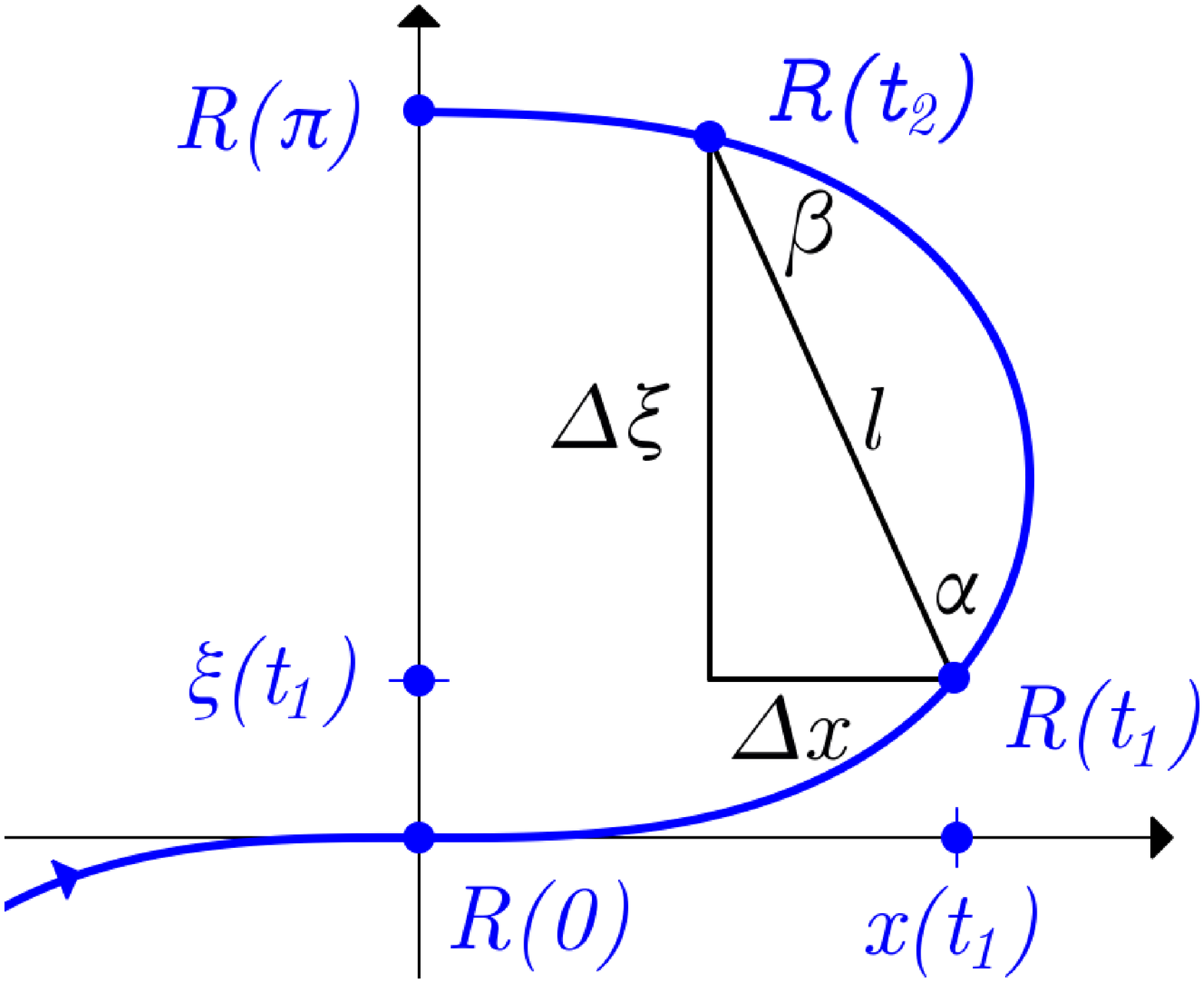}
\newline
Assuming $t_1<t_2$, we  introduce the following
notation (see Fig.\ \figseven):
$$
\Delta x:=\sin (t_2)-\sin (t_1),\quad \Delta \xi:=\xi (t_2)-\xi(t_1),\quad l:=|R(t_2)-R(t_1)|,
$$
whereby $l^2=(\Delta x)^2+(\Delta \xi )^2$ and $\norm{R_{[t_1,t_2]}}^2=\Delta \xi$.
We refer to the quantity $l\norm{R_{[t_1,t_2]}}^2$ as the 
{\bf normalized bending energy} of $R_{[t_1,t_2]}$ because this would be the resultant
bending energy if  $R_{[t_1,t_2]}$ were scaled by the factor $1/l$. Note that if  $R_{[t_1,t_2]}$
is similar to a curve in $\Sopt(\ga,\gb)$ (defined in Section 3), then we have
$$
E(\ga,\gb)=l\norm{R_{[t_1,t_2]}}^2 = l \Delta \xi.
$$

Let $Q$ denote
the mapping $(t_1,t_2)\mapsto (\ga,\gb)$ so that
$$
(\ga,\gb)=Q(t_1,t_2).
$$ 
We leave it to the reader to verify the following formulae for partial derivatives (these are valid for any
sufficiently smooth curve):
$$
\aligned
{\partial \ga\over \partial t_1}&=|R'(t_1)|\bigg({\sin \ga\over l}+\kappa (t_1)\bigg) \\
{\partial \gb\over\partial t_1}&=|R'(t_1)|{\sin \ga \over l} 
\endaligned\qquad
\aligned
{\partial \ga\over \partial t_2}&=-|R'(t_2)|{\sin \gb \over l} \\
{\partial \gb\over \partial t_2}&=|R'(t_2)|\bigg({-\sin \gb \over l}+\kappa (t_2)\bigg)
\endaligned
\tag\Elabel\elalphar$$
The determinant of $DQ:= \bmatrix {\partial \ga\over \partial t_1}&{\partial
\ga\over \partial t_2}\\{\partial \gb\over \partial
t_1}&{\partial \gb\over \partial t_2}\endbmatrix$ is therefore given by
$$
\det(DQ)=|R'(t_1)||R'(t_2)|\pr{\kappa (t_1)\kappa (t_2)+\kappa
(t_2){\sin \ga \over l}-\kappa (t_1){\sin \gb \over l}}.
\tag\Elabel\detDQgeneric$$
Let the {\it cross product} in $\BbC$ be defined by $(u_1+iv_1)\times (u_2+iv_2):=u_1v_2-v_1u_2$. Noting that
$l |R'(t_1)|\sin\ga = (R(t_2)-R(t_1))\times R'(t_1) = -\cos t_1\gD\xi +\xi'(t_1)\gD x $ and 
$l|R'(t_2)|\sin\gb = (R(t_2)-R(t_1))\times R'(t_2) = -\cos t_2\gD\xi + \xi'(t_2)\gD x$, 
the generic formulation in (\detDQgeneric) can be written specifically as:
$$
\aligned \det (DQ)&={4\sin t_1\sin t_2\over \sqrt {1+\sin^2t_1}\sqrt{1+\sin^2t_2}}
+ {2 \sin t_2\over l^2 \sqrt {1+\sin^2 t_2}}\big( -\cos t_1\,\Delta \xi + \xi'(t_1)\Delta x\big)\\
&-{2 \sin t_1\over l^2\sqrt {1+\sin ^2 t_1}} \big(- \cos t_2\,\Delta \xi + \xi'(t_2)\Delta x\big).
\endaligned
\tag\Elabel\detDQspecial$$
Note that if both $\sin t_1=0$ and $\sin t_2=0$, then $\det(DQ)=0$.
\proclaim{Lemma \Tlabel\bdrycase} Suppose $(t_1,t_2)$ belongs to the first or third set defined in
Theorem \detq.  If $\sin t_1 \sin t_2=0$, then $\det(DQ)<0$.
\endproclaim
\demo{Proof} We prove the lemma assuming $t_1=0<t_2<\pi$ since the proof in the other three cases is similar.
Since $\x'(0)=0$ and $\gD\xi>0$, it follows from (\detDQspecial) that 
$\dsize\det(DQ)={2 \sin t_2\over l^2 \sqrt {1+\sin^2 t_2}}(-\gD\xi)<0$.
\qed\enddemo
If $\sin t_1\sin t_2\neq 0$, then (\detDQspecial) can be factored as
$$
\aligned \det(DQ)&={2\gD\xi\over l^2\sqrt {1+\sin^2t_1}\sqrt
{1+\sin^2t_2}} \sin t_1\sin t_2\, W(t_1,t_2),\quad\text{ where }\\
&W(t_1,t_2):=2\Delta \xi +{(\Delta x)^2\over \Delta
\xi}+{\cos t_2\sqrt {1+\sin ^2 t_2}\over \sin t_2}-{\cos t_1\sqrt
{1+\sin ^2 t_1}\over \sin t_1}.
\endaligned
\tag\Elabel\eldefg
$$
Note that the sign of $\det(DQ)$ is the same as that of $\sin t_1\sin t_2\, W(t_1,t_2)$.
\proclaim{Lemma \Tlabel\partialG} If $\sin t_1\sin t_2\neq 0$, then
$$\align
{\partial W\over \partial t_1}&={\sqrt {1+\sin^2t_1}\over (\Delta
\xi)^2}\bigg[{\cos t_1\Delta \xi\over \sin t_1}-{\sin t_1\Delta
x\over \sqrt {1+\sin^2t_1}}\bigg]^2\geq 0,\quad\text{ and}\\
{\partial W\over \partial t_2}&=-{\sqrt {1+\sin^2t_2}\over (\Delta
\xi)^2}\bigg[{\cos t_2\Delta \xi\over \sin t_2}-{\sin t_2\Delta
x\over \sqrt {1+\sin^2t_2}}\bigg]^2\leq 0.
\endalign
$$
\endproclaim
\demo{Proof} We only prove the result pertaining to ${\partial W\over \partial t_1}$ since
the proof of the other is the same, {\it mutatis mutandis}. Direct differentiation yields
$$\eqalign{ {\partial W\over \partial t_1}=-2\x'(t_1)
+&{-2\Delta x\Delta \xi\cos t_1 +(\Delta x)^2\x'(t_1)\over
(\Delta \xi)^2}\cr \noalign {\vskip 5pt} -&{-\sin ^2t_1\sqrt
{1+\sin ^2 t_1}+{\cos ^2t_1\sin^2t_1\over \sqrt
{1+\sin^2t_1}}-\cos^2t_1\sqrt{1+\sin^2t_1}\over \sin^2t_1 },\cr }$$
which then simplifies to
$$\eqalign{{\partial W\over \partial t_1}&={-2\cos t_1\Delta
x\Delta \xi +{\sin ^2t_1\over \sqrt {1+\sin ^2 t_1}}(\Delta
x)^2\over (\Delta \xi)^2}+\bigg({\sqrt {1+\sin^2 t_1}\over
\sin^2t_1}-{1+\sin^2t_1\over \sqrt {1+\sin^2 t_1}}\bigg)\cr
\noalign {\vskip 5pt} &={-2\cos t_1\Delta x\Delta \xi +{\sin
^2t_1\over \sqrt {1+\sin ^2 t_1}}(\Delta x)^2\over (\Delta
\xi)^2}+{\cos^2t_1\sqrt {1+\sin^2t_1}\over \sin^2t_1}.\cr}$$
A simple computation then shows that this last expression can be
factored as stated in the lemma.
\qed\enddemo
\proclaim {Theorem \Tlabel \detq} There exists a unique $t^*\in
(0,\pi)$ such that $W(-t^*,t^*)=0$.
Moreover $\det (DQ)<0$ on the following sets:\newline
(i) $\set{(t_1,t_2):-\pi\leq t_1<t_2\leq 0,\;(t_1,t_2)\neq (-\pi,0)}$,\newline
(ii) $\set{(t_1,t_2):-t^*<t_1<0<t_2<t^*}$\newline
(iii) $\set{(t_1,t_2):0\leq t_1<t_2\leq \pi,\;(t_1,t_2)\neq (0,\pi)}$,\newline
(iv) $\set{(t_1,t_2):\pi-t^*<t_1<\pi<t_2<\pi+t^*}$\newline
\endproclaim
\demo {Proof} 
For $-\pi<t_1<0<t_2<\pi$, the function $W(t_1,t_2)$ is analytic in both $t_1$ and $t_2$,
and consequently, it follows from Lemma \partialG\ that $W(t_1,t_2)$ is  
increasing in $t_1$ and   decreasing in $t_2$. Furthermore, the function
$W(-t,t)$ is analytic and   decreasing for $t\in(0,\pi)$. Note that if
$-\frac\pi2\leq t_1<0<t_2\leq \frac\pi2$, then $\sin t_1<0$ and it is clear
(from (\eldefg)) that $W(t_1,t_2)>0$. In particular, $W(-t,t)>0$ for all $t\in(0,\frac\pi2]$.
It is easy to verify
(by inspection of (\eldefg)) that $\lim_{t\to \pi ^-}W(-t,t)=-\infty$, and so it follows that
there exists a unique $t^*\in(0,\pi)$ such that $W(-t^*,t^*)=0$. 

If $(t_1,t_2)$ belongs to set (ii), then $W(t_1,t_2)>W(-t^*,t_2)>W(-t^*,t^*)=0$ and
since $\sin t_1 \sin t_2<0$, it follows that $\det(DQ)<0$. This proves that
$\det(DQ)<0$ for all $(t_1,t_2)$ in set (ii).

We will  show that $\det(DQ)<0$ for all $(t_1,t_2)$  in set (i).
This has already been proved in Lemma \bdrycase\ if $0=t_1<t_2<\pi$ or
$0<t_1<t_2=\pi$, so assume $0<t_1<t_2<\pi$. As above,
the function $W(t,t_2)$ is analytic and   increasing for $t\in(0,t_2)$.
It is easy to see (by inspection of (\eldefg)) that 
$\lim _{t\to t_2^-}W(t,t_2)=0$, and
therefore $W(t,t_2)<0$ for all $t\in(0,t_2)$; in particular, $W(t_1,t_2)<0$.
Since $\sin t_1 \sin t_2>0$, we have $\det(DQ)<0$.
This completes the proof that $\det(DQ)<0$ for all $(t_1,t_2)$ in set (i).

Finally, if $(t_1,t_2)$ belongs to set (iii) or set (iv), then 
$(t_1-\pi,t_2-\pi)$ belongs to set (i) or set (ii) and
$\det(DQ(t_1,t_2))=\det(DQ(t_1-\pi,t_2-\pi))<0$.
\qed\enddemo
\proclaim {Corollary \Tlabel \cortstar} Let $t^*\in(0,\pi)$ be as defined in Theorem \detq.
Then $t^*>\frac\pi2$ and $\gb(0,t^*)>\frac\pi2$. Moreover, $\gb(0,t)$ is 
  increasing for $t\in(0,t^*]$ and   decreasing 
for $t\in[t^*,\pi]$.
\endproclaim
\demo{Proof}
Since $W(-t,t)>0$ for $t\in(0,\frac\pi2]$, it is clear that $t^*>\frac\pi2$.
Since $W(-t^*,t^*)=0$, it follows from (\eldefg) that  $\det(DQ(-t^*,t^*))=0$, 
and therefore, by (\elalphar), we must have 
$$
\kappa (-t^*)\kappa (t^*)+\kappa
(t^*){\sin (\ga(-t^*,t^*)) \over l(-t^*,t^*)}-\kappa (-t^*){\sin (\gb(-t^*,t^*)) \over l(-t^*,t^*)}=0.
$$
From the definition of $\ga$ and $\gb$ it is clear that
$\ga (-t^*,t^*)=\gb (-t^*,t^*)>0$ and $\kappa (t^*)=-\kappa
(-t^*)>0$, so the above equality reduces to
$\dsize\kappa (t^*)-{2\sin (\gb (-t^*,t^*))\over l(-t^*,t^*)}=0$.
From the symmetry of the curve $R$ one has $\sin (\gb
(-t^*,t^*))=\sin (\gb (0,t^*))$ and $l(-t^*,t^*)=2l(0,t^*)$ which yields
$\kappa (t^*)-{\sin (\gb (0,t^*))\over l(0,t^*)}=0$.
It now follows from (\elalphar) that
$\dsize{\partial \gb \over \partial t_2}(0,t^*)=0$. Moreover, the
uniqueness of $t^*\in (0,\pi)$ shows (running the above argument
backwards) that $t=t^*$ is the unique $t\in (0,\pi)$ where
$\dsize{\partial \gb \over \partial t_2}(0,t)=0$. This implies that
the function $\gb (0,t)$ is   increasing on $(0,t^*]$ and
decreasing on $[t^*,\pi]$. Consequently,  $\gb (0,t^*)>\gb(0,\pi)=\frac\pi2$.
\qed\enddemo
\proclaim{Corollary \Tlabel\cortbar}
There exists a unique $\tbar\in(0,t^*)$ such that  $\gb(0,\tbar)=\frac\pi2$.
Moreover, we have $\gb(0,t)<\frac\pi2$ for all $0<t<\tbar$ and
$\gb(0,t)>\frac\pi2$ for all $\tbar<t<\pi$.
\endproclaim
\demo{Proof}
Since $\lim_{t\to 0^+}\gb(0,t)=0$, $\gb(0,t^*)> \frac\pi2$, and  $\gb(0,\pi)=\frac\pi2$, the result
follows immediately from Corollary \cortstar.
\qed\enddemo
%
%
\Section{Unicity of Parameters}
%
For $(\ga,\gb)\in[-\frac\pi2,\frac\pi2]^2$, recall that $c_1(\ga,\gb)$ is the unique $C^\infty$ s-curve in $\Sopt(\ga,\gb)$.
In Theorem \sunique\ (iii), it is shown that if $(\ga,\gb)\neq (0,0)$, then there exist $t_1<t_2<t_1+2\pi$ such
that $c_1(\ga,\gb)$ is directly similar to $R_{[t_1,t_2]}$.  In this section, we are concerned with the unicity of
the parameters $(t_1,t_2)$. 
The rectangular elastic curve $R$ is periodic in the sense that $R(t+2\pi)=i2d+R(t)$, and it follows
that $R_{[t_1',t_2']}$ is directly congruent to $R_{[t_1,t_2]}$ whenever $(t_1',t_2')=(t_1,t_2)+k(2\pi,2\pi)$ for some
integer $k$; in particular $Q(t_1',t_2')=Q(t_1,t_2)$.  With the identification $(t_1',t_2')\equiv(t_1,t_2)$, the half-plane 
$Y:=\set{(t_1,t_2):t_1\leq t_2}$ becomes a half-cylinder, with boundary $t_1=t_2$, and we adopt the view that $Q$ is defined
on the interior of the cylinder $Y$. 

In this section, we will prove the following.
\proclaim{Theorem \Tlabel\thuniqueone} For all $(\ga,\gb)\in[-\frac\pi2,\frac\pi2]^2\bs\set{(0,0)}$, there exists a unique $(t_1,t_2)$
in the cylinder $Y$ such that $t_1<t_2<t_1+2\pi$ and $R_{[t_1,t_2]}$ is an s-curve with chord angles $(\ga,\gb)$.
\endproclaim
\proclaim{Theorem \Tlabel\thuniquetwo} Let $t_1<t_2<t_1+2\pi$ be such that $R_{[t_1,t_2]}$ is an s-curve with chord angles
$(\ga,\gb)\in[-\frac\pi2,\frac\pi2]^2$. Then $R_{[t_1,t_2]}$ is directly similar to $c_1(\ga,\gb)$.
\endproclaim
%
We define the following subsets of the interior of $Y$:
$$
\matrix U_0:=\set{(t_1,t_2): -\pi\leq t_1<t_2\leq 0} & V_1:=\set{(t_1,t_2): -\pi\leq t_1<0<t_2\leq \pi} \\
U_2:=\set{(t_1,t_2):0\leq t_1<t_2\leq\pi} & V_3:=\set{(t_1,t_2): 0\leq t_1<\pi<t_2\leq 2\pi} \endmatrix.
$$

These sets are pairwise disjoint subsets of the cylinder $Y$, and
for $t_1<t_2$, it is easy to verify that $R_{[t_1,t_2]}$ is a right c-curve if and only if
$(t_1,t_2)\in U_0$,
a non-degenerate right-left s-curve if and only if $(t_1,t_2)\in V_1$,
a left c-curve if and only if $(t_1,t_2)\in U_2$, and
a non-degenerate left-right s-curve if and only if $(t_1,t_2)\in V_3$. 

The restriction $t_1<t_2<t_1+2\pi$ eliminates $(-\pi,\pi)$ from $V_1$ and $(0,2\pi)$ from $V_3$, and we therefore have 
the following as a consequence of Theorem \sunique\ (iii).
\proclaim{Proposition \Tlabel\anticipation} For all $(\ga,\gb)\in[-\frac\pi2,\frac\pi2]^2\bs\set{(0,0)}$, there exists \newline
$(t_1,t_2)\in U_0\cup V_1\cup U_2\cup V_3\bs\set{(-\pi,\pi),(0,2\pi)}$ such that $c_1(\ga,\gb)$ is directly similar to
$R_{[t_1,t_2]}$.
\endproclaim
In particular, we have the following corollary.
%
\proclaim{Corollary \Tlabel\ontoone} For all $(\ga,\gb)\in[-\frac\pi2,\frac\pi2]^2\bs\set{(0,0)}$, there exists \newline
$(t_1,t_2)\in U_0\cup V_1\cup U_2\cup V_3\bs\set{(-\pi,\pi),(0,2\pi)}$ such that $Q(t_1,t_2)=(\ga,\gb)$.
\endproclaim
We intend to show that the pair $(t_1,t_2)$ is unique, but before beginning the proof of this, we will harmlessly replace
$V_1,V_3$ with smaller sets $U_1,U_3$, defined below.

With $\tbar$ as defined in Corollary \cortbar, we define
$$\matrix U_1:=\set{(t_1,t_2): -\tbar\leq t_1<0<t_2\leq \tbar} & U_3:=\set{(t_1,t_2): \pi-\tbar\leq t_1<\pi<t_2\leq \pi+\tbar}
\endmatrix. $$
%
\proclaim{Lemma \Tlabel\notorius} If $(t_1,t_2)$ belongs to $V_1\bs U_1$ or $V_3\bs U_3$ and satisfies $t_2-t_1<2\pi$, then
$(\ga,\gb)\not\in [-\frac\pi2,\frac\pi2]^2$.
\endproclaim
\demo{Proof} We will only prove the lemma for $V_1\bs U_1$ since the proof for $V_3\bs U_3$ is similar.
Let $(t_1,t_2)\in V_1\bs U_1$ satisfy $t_2-t_1<2\pi$. We can assume, without loss of generality, that $t_2\geq -t_1$, since
the remaining case $t_2<-t_1$ is similar. We will show that $\gb>\frac\pi2$. If $t_2=-t_1$, then we must have $\tbar<t_2<\pi$
and, by symmetry, $\gb=\gb(0,t_2)$;  hence $\gb=\gb(0,t_2)>\frac\pi2$ by Corollary \cortbar.  
So assume $t_2>-t_1$, which implies $\tbar<t_2\leq \pi$. The chord $[R(t_1),R(t_2)]$
must intersect the negative $x$-axis, since otherwise we would have $t_2\leq -t_1$.  Therefore, $\gb>\gb(0,t_2)>\frac\pi2$.
\qed \enddemo
As a consequence of the lemma, the set $U_0\cup V_1\cup U_2\cup V_3$ in Corollary \ontoone\ can be replaced with
$U:=U_0\cup U_1\cup U_2\cup U_3$:
\proclaim{Corollary \Tlabel\quickrain} For all $(\ga,\gb)\in[-\frac\pi2,\frac\pi2]^2\bs\set{(0,0)}$,
there exists $(t_1,t_2)\in U$ such that  $Q(t_1,t_2)=(\ga,\gb)$.
\endproclaim
\noindent
\epsfxsize=2truecm\epsffile{space.eps}\epsfysize=6truecm \epsffile{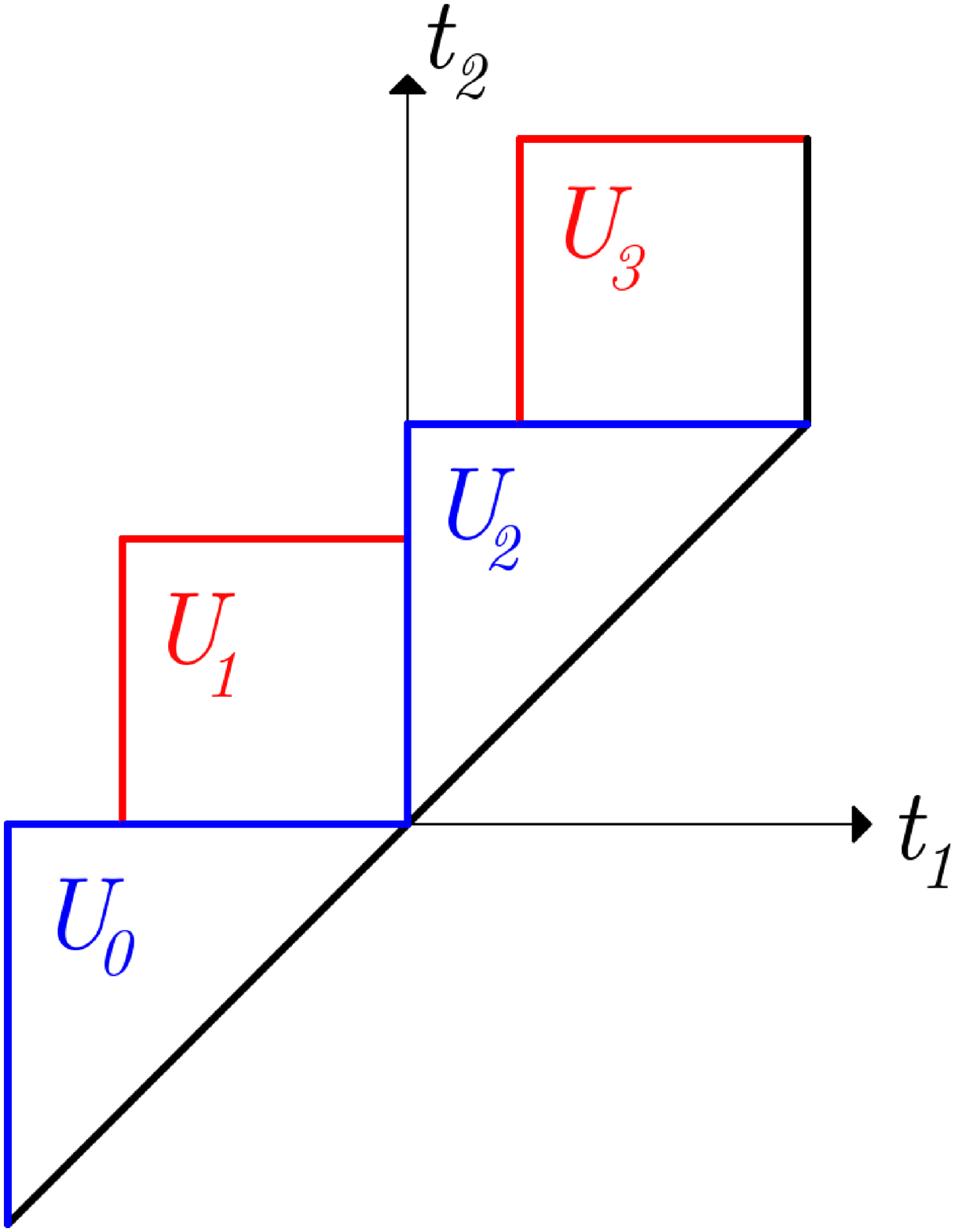} \hskip2truecm \epsfysize=6truecm \epsffile{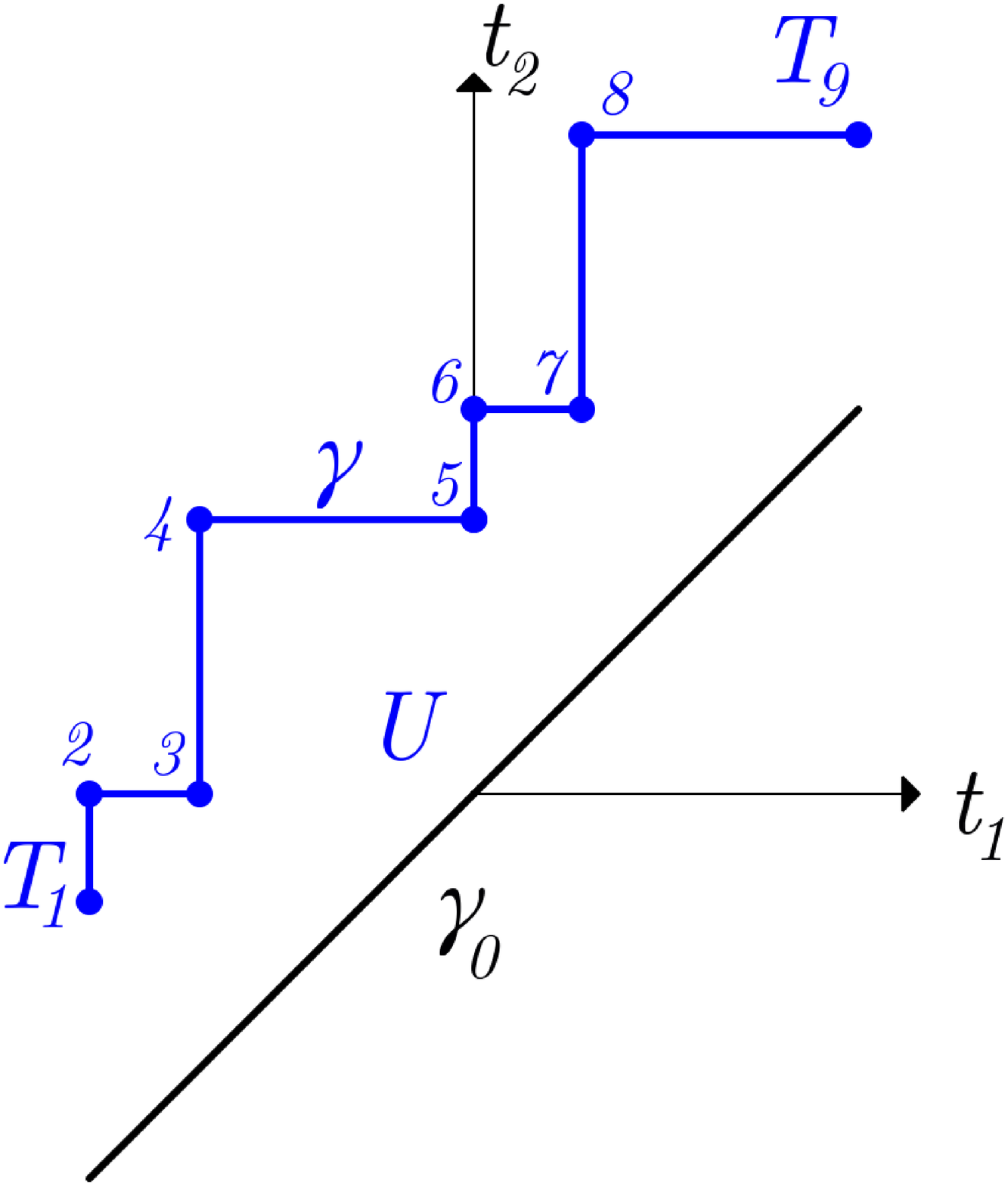} 
\newline
{\bf Fig.\ \Flabel\figten} \ \ \ {\bf (a)}  the sets  $U_0, U_1, U_2, U_3$\ \ \ \ \ \ \  \ \ \ \ \ \ \ \  
{\bf (b)} the set $U$ and its boundary $\gga_0\cup \gga$

In Fig.\ \figten (a), the sets  $U_0, U_1, U_2, U_3$ are depicted on the fundamental domain
$-\pi \leq t_1 <\pi$ of the cylinder $Y$, and their union $U$ is depicted in Fig.\ \figten(b).  The set $U$ is bounded
below by the line $\gga_0:=\set{(t_1,t_2):t_1=t_2}$ (which is not contained in $U$) and above by the staircase path
$\gga:=[T_1,T_2,\ldots,T_9]$ (which is contained in $U$). Here,
$T_1=(-\pi,\tbar-\pi)$, $T_2=(-\pi,0)$, $T_3=(-\tbar,0)$, $T_4=(-\tbar,\tbar)$, and
$T_i=T_{i-4}+(\pi,\pi)$ for $i=5,6,7,8,9$.  Note that on the cylinder $Y$,
the vertical half line starting from $\gga_0$ and passing through $T_9$ is identified with the same, but passing through $T_1$;
in particular $T_9$ is identified with $T_1$.

At present, $Q$ is defined and is $C^\infty$ on the interior of the cylinder $Y$. On the boundary of $Y$ (the line $\gga_0$),
we define $Q$ to be $(0,0)$; in other words, we define $\ga(t,t):=0$ and $\gb(t,t):=0$ for all $t\in\BbR$.
\proclaim{Lemma \Tlabel\december} $Q$ is continuous on the cylinder $Y$.
\endproclaim
\demo{Proof} We will show that
$\abs{\ga(t_1,t_2)}+\abs{\gb(t_1,t_2)}\leq 2(t_2-t_1)$ whenever $t_1<t_2$. It is generally true that the absolute
sum of the chord angles is bounded by the absolute turning angle of the curve. In the present context, this means
that $\abs\ga + \abs\gb \leq \int_{t_1}^{t_2} \abs{\gk(t)}\,\abs{R'(t)}\,dt$. Since $\abs{\gk(t)}=\abs{2\sin t}\leq 2$
and $\abs{R'(t)}=1/\sqrt{1+\sin^2 t}\leq 1$, the desired inequality is immediate.
\qed\enddemo
\noindent
\epsfxsize=1truecm\epsffile{space.eps}\epsfysize=6truecm \epsffile{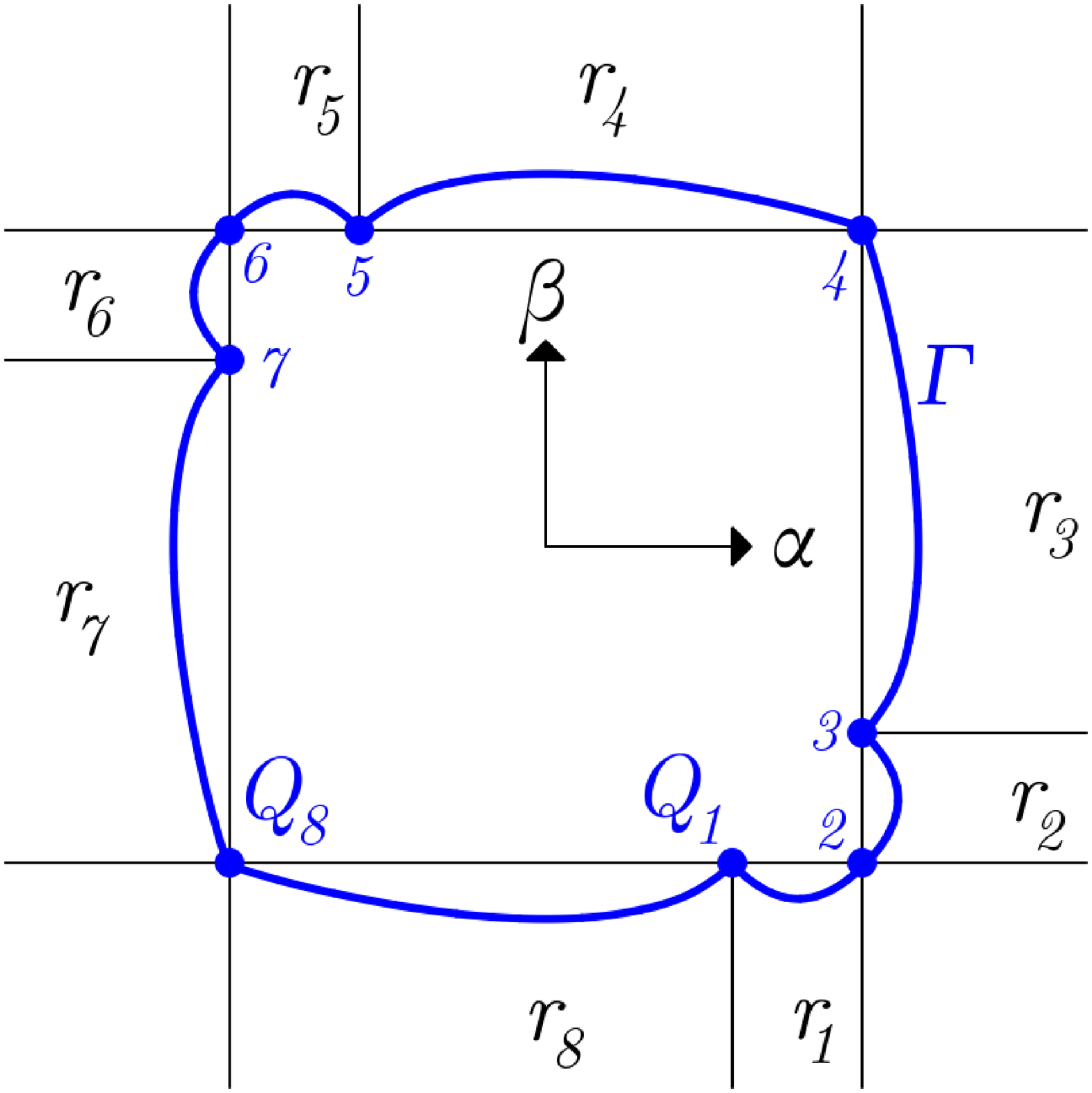} \hskip2truecm \epsfysize=6truecm \epsffile{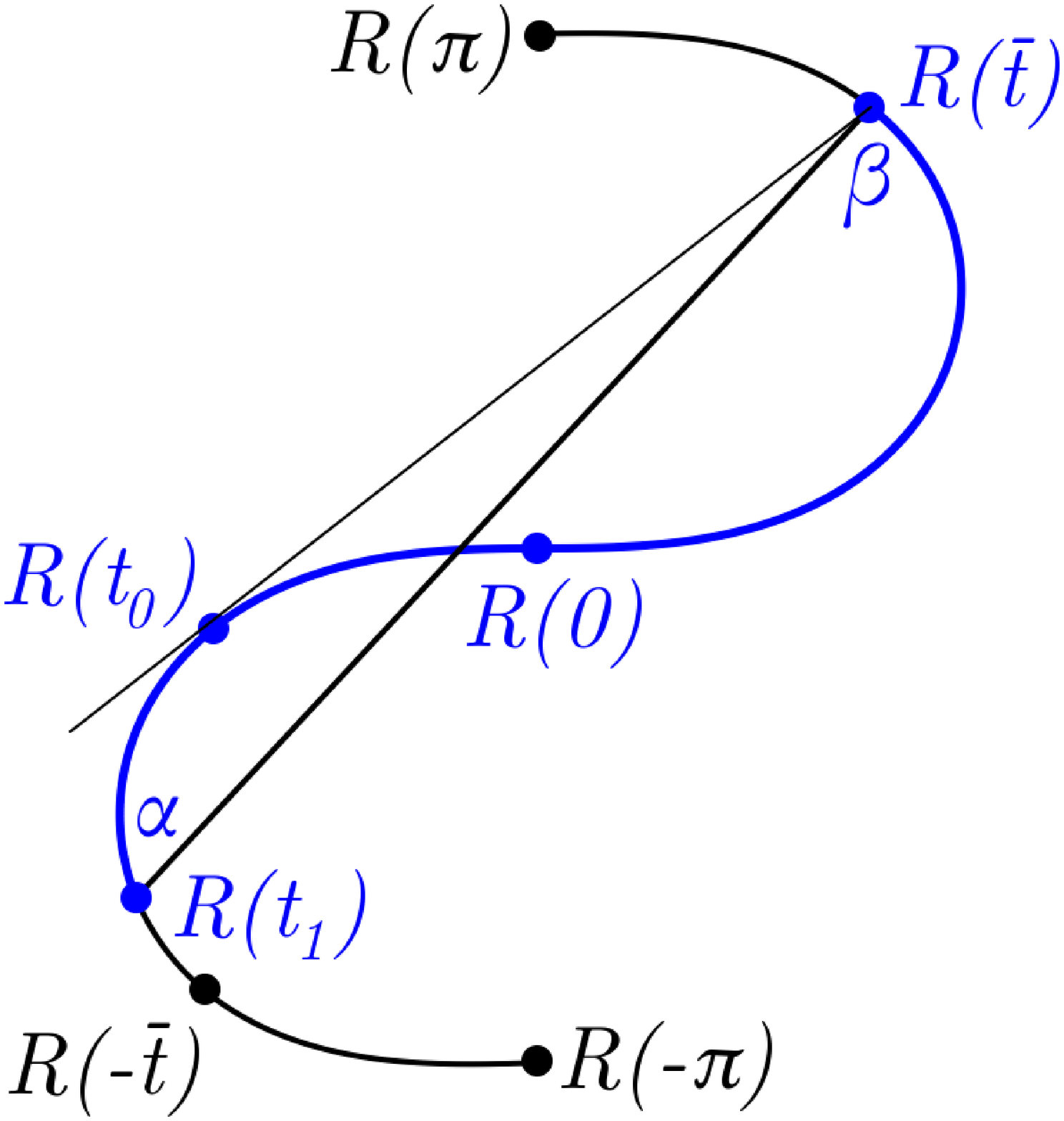}
\newline
{\bf Fig.\ \Flabel\figeight} \ \ \ \ {\bf (a)} the image $\gG:=Q(\gga)$ \ \ \ \ \ \ \ \  \ \ \ \ \ \ \
{\bf (b)} the parameters $-\tbar<t_1<t_0<0$

Fig.\ \figeight\ (a) depicts the image  $\gG:=Q(\gga)$ where $Q_i:=Q(T_i)$ are given by
$Q_1=(\psibar,-\frac\pi2)$, $Q_2=(\frac\pi2,-\frac\pi2)$, $Q_3=(\frac\pi2,-\psibar)$, $Q_4=(\frac\pi2,\frac\pi2)$, and
$Q_i=-Q_{i-4}$ for $i=5,6,7,8,9$; here $\psibar:=|\ga(0,\tbar)|$.
The staircase path $\gga$ consists of eight segments
$[T_i,T_{i+1}]$, $i=1,2,\ldots,8$, and it is apparent in Fig.\ \figeight(a) that their images $\set{Q([T_i,T_{i+1}])}$
belong to eight non-overlapping unbounded rectangles $\set{r_i}$. Specifically,
$r_1:=[\psibar,\frac\pi2]\times (-\infty,-\frac\pi2]$, $r_2:=[\frac\pi2,\infty)\times [-\frac\pi2,-\psibar]$,
$r_3:=[\frac\pi2,\infty) \times [-\psibar,\frac\pi2]$, $r_4:=[-\psibar,\frac\pi2]\times [\frac\pi2,\infty)$, and
$r_i:=-r_{i-4}$ for $i=5,6,7,8$.
\proclaim{Lemma} For $i=1,2,\ldots,8$,
$Q$ is injective on $[T_i,T_{i+1}]$ and maps the interior of $[T_i,T_{i+1}]$ into the interior of $r_i$.
\endproclaim
\demo{Proof} Let us first consider the case $i=4$, where $r_4=[-\psibar,\frac\pi2]\times [\frac\pi2,\infty)$.
Along the segment $[T_4,T_5]$ (see Fig.\ \figten(b)), $t_1$ ranges from $-\tbar$ to $0$, while
$t_2=\tbar$ is fixed. At the endpoints, we have $\ga(-\tbar,\tbar)=\gb(-\tbar,\tbar)=\frac\pi2$ and
$\ga(0,\tbar)=-\psibar$, $\gb(0,\tbar)=\frac\pi2$.
Since $\gk(t)<0$ for $t\in(-\pi,0)$, it is clear (see Fig.\ \figeight(b)) that
$\gb(t_1,\tbar)>\gb(0,\tbar)=\frac\pi2$ and $\ga(t_1,\tbar)<\ga(t_1,0)=\gb(0,-t_1)<\frac\pi2$ for $t_1\in(-\tbar,0)$.\newline
Recall from (\elalphar) that $\frac{\partial\ga}{\partial t_1}=|R'(t_1)|\pr{\frac{\sin\ga}\ell + \gk(t_1)}$.
Note that if  $t_1\in(-\tbar,0)$ and $\ga(t_1,\tbar)\leq 0$, then
$\frac{\partial\ga}{\partial t_1}<0$. From this, one easily deduces that there exists $t_0\in(-\tbar,0)$ such that
$\ga(t_1,\tbar)>0$ for $t_1\in(-\tbar,t_0)$ and $\ga(t_1,\tbar)<0$ for $t_1\in(t_0,0]$. Furthermore,
$\ga(t_1,\tbar)$ is decreasing for $t_1\in[t_0,0]$, and therefore we have
$\ga(t_1,\tbar)\in (-\psibar,\frac\pi2)$ for $t_1\in(-\tbar,0)$.  This completes the proof that $Q$ maps the interior of
$[T_4,T_5]$ into the interior of $r_4$. We will now show that $Q$ is injective on $[T_4,T_5]$. Recall from (\elalphar)
that $\frac{\partial\gb}{\partial t_1}=|R'(t_1)|\frac{\sin\ga}\ell$ and hence
$\gb(t_1,\tbar)$ is increasing when $\ga$ is positive (i.e., for $t_1\in[-\tbar,t_0)$) and
$\gb(t_1,\tbar)$ is decreasing when $\ga$ is negative (i.e., for $t_1\in(t_0,0]$).  Consequently,
$Q$ is injective on $[T_4,T_5]$. This proves the lemma in the case $i=4$ and the cases $i=3,7,8$ follow
by symmetry.\newline
We next consider the case $i=5$, where $r_5=[-\frac\pi2,-\psibar]\times [\frac\pi2,\infty)$. Along the segment $[T_5,T_6]$
(see Fig.\ \figten(b)), $t_1=0$ is fixed while $t_2$ ranges from $\tbar$ to $\pi$.
It is shown in Corollary \cortbar\ that
$\gb(0,t_2)>\frac\pi2$ for all $t_2\in(\tbar,\pi)$. Recall from (\elalphar)
that  $\frac{\partial\ga}{\partial t_2}=-|R'(t_2)|\frac{\sin\gb}\ell$. Since $\gb(0,t_2)>0$, it follows that
$\frac{\partial\ga}{\partial t_2}<0$ for all $t_2\in[\tbar,\pi]$ and hence $\ga(0,t_2)$ is decreasing for
$t_2\in[\tbar,\pi]$. Consequently, $Q$ is injective on $[T_5,T_6]$, and since $\ga(0,\tbar)=-\psibar$ and
$\ga(0,\pi)=-\frac\pi2$, it also follows that $Q$ maps the interior of $[T_5,T_6]$ into the interior of $r_5$.
This proves the lemma for the case $i=5$ and the cases $1,2,6$ follow by symmetry.
\qed\enddemo
\proclaim{Proposition \Tlabel\injectgamma} The following hold.\newline
(i) $Q$ is continuous on $U\cup \gga_0$.\newline
(ii) In the interior of $U$, $Q$ is $C^\infty$ and its Jacobian is nonzero.\newline
(iii) $Q(\gga_0)=\set{(0,0)}$ and $Q(t_1,t_2)\neq (0,0)$ for all $(t_1,t_2)\in U$\newline
(iv) $Q$ is injective on $\gga$.
\endproclaim
\demo{Proof} Item (i) is a consequence of Lemma \december, and (ii) is proved in Theorem \detq.
The first assertion in (iii), $Q(\gga_0)=\set{(0,0)}$, holds by definition.  It is easy to verify that if $f$ is
an s-curve with chord angles $(\ga,\gb)=(0,0)$, then $f$ is a line segment. But $R_{[t_1,t_2]}$ is
never a line segment because the signed curvature of $R$ only vanishes at times $k\pi$, $k\in\BbZ$. Since
$R_{[t_1,t_2]}$ is an s-curve for all $(t_1,t_2)\in U$, we obtain the second assertion in (iii).  Since the rectangles
$r_1,r_2,\ldots,r_8$ are non-overlapping, we obtain (iv) as a consequence of the above lemma.
\qed\enddemo
On the basis of Proposition \injectgamma, we have the following, which is proved in the Appendix.
\proclaim{Theorem \Tlabel\inject} $Q$ is injective on $U$.
\endproclaim
\remark{Remark} The proof of Theorem \inject\ can be extended to show that $Q$ is injective on the larger set $U$ obtained
with $U_1$ and $U_3$ defined with $t^*$ in place of $\tbar$.
\endremark
We can now easily prove Theorems \thuniqueone\ and \thuniquetwo.
\demo{Proof of Theorem \thuniqueone}
Let $(\ga,\gb)\in[-\frac\pi2,\frac\pi2]^2\bs\set{(0,0)}$.  It follows from Corollary \quickrain\ and Theorem \inject\ that
there exists a unique $(t_1,t_2)\in U$ such that $Q(t_1,t_2)=(\ga,\gb)$; this establishes existence.  Now, if
$(t_1,t_2)\in Y$ is such that $t_1<t_2<t_1+2\pi$ and $R_{[t_1,t_2]}$ is an s-curve with chord angles $(\ga,\gb)$, then
it follows from Lemma \notorius\ and the observations made above Proposition \anticipation\ that $(t_1,t_2)\in U$,
whence follows uniqueness.
\qed\enddemo
\demo{Proof of Theorem \thuniquetwo} Assume $t_1<t_2<t_1+2\pi$ and
that $R_{[t_1,t_2]}$ is an s-curve. From the observations above
Proposition \anticipation, it follows that $(t_1,t_2)$, as a point
on the cylinder $Y$, belongs to $U_0\cup V_1 \cup U_2 \cup
V_3\bs\set{(-\pi,\pi),(0,2\pi)}$. Assume that the chord angles
$(\ga,\gb)$ of $R_{[t_1,t_2]}$ belong to $[-\frac\pi2,\frac\pi2]$.  As
mentioned in the proof of Proposition \injectgamma\ (iii), we must
have $(\ga,\gb)\neq(0,0)$ and therefore, by Proposition
\anticipation, there exists $(t_1',t_2')\in U_0\cup V_1 \cup U_2
\cup V_3\bs\set{(-\pi,\pi),(0,2\pi)}$ such that $c_1(\ga,\gb)$ is
directly similar to $R_{[t_1',t_2']}$. Since
$Q(t_1,t_2)=(\ga,\gb)=Q(t_1',t_2')$, it follows from Theorem
\thuniqueone\ that $(t_1,t_2)=(t_1',t_2')$ (in the cylinder $Y$)
and therefore $R_{[t_1',t_2']}$ is directly congruent to
$R_{[t_1,t_2]}$; hence $R_{[t_1,t_2]}$ is directly similar to
$c_1(\ga,\gb)$. \qed\enddemo
%
\Section {Proof of Condition (\consistent) }   
%
In this section we prove that condition (\consistent) holds with $\mu=2$:
\proclaim{Theorem \Tlabel\Tconsistent} For  all 
$(\ga_0,\gb_0)\in[-\frac\pi2,\frac\pi2]^2\bs\set{(-\frac\pi2,\frac\pi2),(\frac\pi2,-\frac\pi2)}$,
$$
[-\gk_a(c_1(\ga_0,\gb_0)), \gk_b(c_1(\ga_0,\gb_0))]=2\nabla E_1(\ga_0,\gb_0).
\tag\Elabel\nablaEQ$$
\endproclaim
\demo{Proof} 
Fix $(\ga_0,\gb_0)\in [-\frac\pi2,\frac\pi2]^2\bs\set{(-\frac\pi2,\frac\pi2),(\frac\pi2,-\frac\pi2)}$. 
We first address the easy case $(\ga_0,\gb_0)=(0,0)$, where $c_1(0,0)$ is a line segment. In the proof
of \cite{\BJone, Prop. 7.6}, it is shown that there exists a constant $C$ such that
$E_1(\ga,\gb)=E(\ga,\gb)\leq C (\tan^2\ga + \tan\ga \tan\gb + \tan^2 \gb)$ for all $(\ga,\gb)\in[-\pi/3,\pi/3]^2$.
From this it easily follows that $\nabla E_1(0,0)=[0,0]$, and since the line segment $c_1(0,0)$ has
$0$ curvature, we obtain (\nablaEQ) for the case $(\ga_0,\gb_0)=(0,0)$. \newline
We proceed assuming $(\ga_0,\gb_0)\in [-\frac\pi2,\frac\pi2]^2\bs\set{(-\frac\pi2,\frac\pi2),(\frac\pi2,-\frac\pi2),(0,0)}$.
It follows from Corollary \quickrain\ that there exists $(\gt_1,\gt_2)\in U$ such that
$Q(\gt_1,\gt_2)=(\ga_0,\gb_0)$. The restriction $(\ga_0,\gb_0)\not\in\set{(-\frac\pi2,\frac\pi2),(\frac\pi2,-\frac\pi2)}$ 
ensures that $(\gt_1,\gt_2)\not\in\set{(0,\pi),(-\pi,0)}$, and consequently, it follows from Theorem \detq\ that $DQ(\gt_1,\gt_2)$
is nonsingular. Since $Q$ is $C^\infty$ on the interior of the cylinder $Y$ (defined in Section 6), 
it follows that there exists an open neighborhood $N$ of
$(\gt_1,\gt_2)$ such that $Q$ is injective on $N$, $DQ$ is nonsingular on $N$, $Q(N)$ is an open neighborhood of $(\ga_0,\gb_0)$, 
and $Q^{-1}$ is $C^\infty$ on $Q(N)$. We define $E^*:Q(N)\to[0,\infty)$ as follows.  For $(\ga,\gb)\in Q(N)$,
$$
E^*(\ga,\gb):= l\norm{R_{[t_1,t_2]}}^2,\text{ where }(t_1,t_2):=Q^{-1}(\ga,\gb)\text{ and }l:=\abs{R(t_1)-R(t_2)}.
$$
{\bf Claim.} If $(\ga,\gb)\in Q(N)\cap [-\frac\pi2,\frac\pi2]^2$, then $E^*(\ga,\gb)=E_1(\ga,\gb)$ and
$c_1(\ga,\gb)$ is directly congruent to $\frac1{l}R_{[t_1,t_2]}$.
\demo{proof} Assume $(\ga,\gb)\in Q(N)\cap [-\frac\pi2,\frac\pi2]^2$. Since $Q(t_1,t_2)=(\ga,\gb)$,
it follows from Theorems \thuniqueone\ and \thuniquetwo\ that $c_1(\ga,\gb)$ is directly similar to $R_{[t_1,t_2]}$. Consequently,
$c_1(\ga,\gb)$ is directly congruent to $\frac1l R_{[t_1,t_2]}$ and
$E_1(\ga,\gb) := \norm{c_1(\ga,\gb)}^2 = E^*(\ga,\gb)$, as claimed.
\enddemo
We recall, from Section 2, that the curvature of $R$ is given by $\gk(t)=2\sin t$, and hence
$\gk_a(c_1(\ga,\gb))=2l\sin t_1$ and $\gk_b(c_1(\ga_0,\gb_0))=2l\sin t_2$.
So with the claim in view, in order to establish (\nablaEQ) it suffices to show that
$$
[-l\sin t_1,l\sin t_2]=\nabla E^*(\ga,\gb),\text{ for all }(\ga,\gb)\in Q(N). 
\tag\Elabel\desiredEQ$$
The bending energy of
$R_{[t_1,t_2]}$ (see Section 2) is given by $\norm{R_{[t_1,t_2]})}^2=\xi(t_2)-\xi(t_1)=:\gD\xi$, and hence
$E^*(\ga,\gb)=l \gD\xi$. Defining $\wt E:N\to[0,\infty)$ by $\wt E(t_1,t_2):=l\gD\xi$, we have $\wt E = E^*\circ Q$, and
therefore, since $DQ$ is nonsingular on $N$, (\desiredEQ) is equivalent to
$$
[-l\sin t_1,l\sin t_2]DQ=\nabla \wt E(t_1,t_2),\text{ for all }(t_1,t_2)\in N. 
$$
This can be written explicitly as 
$$\aligned -l\sin t_1\frac{\partial\ga}{\partial t_1}+l\sin t_2\frac{\partial\gb}{\partial t_1} &=\frac{\partial}{\partial t_1}(l\Delta \xi)\\
-l\sin t_1\frac{\partial\ga}{\partial t_2}+l\sin t_2\frac{\partial\gb}{\partial t_2} &=\frac{\partial}{\partial t_2}(l\Delta \xi)
\endaligned\tag\Elabel\eldenergy$$
Using (\elalphar) and the formulae above (\detDQspecial) the first equality is proved as follows.
$$\align
-l\sin t_1\frac{\partial\ga}{\partial t_1}+l\sin t_2\frac{\partial\gb}{\partial t_1} 
&=|R'(t_1)|\sin\ga\gD x - l\sin t_1 |R'(t_1)|\gk(t_1)\\
&=(-\cos t_1 \gD\xi + \xi'(t_1)\gD x)\frac{\gD x}l - 2l\xi'(t_1) \\
&=(-\cos t_1 \gD\xi + \xi'(t_1)\gD x)\frac{\gD x}l - \frac{\gD x^2 + \gD\xi^2}l\xi'(t_1) - l\xi'(t_1)\\
&=-\frac{-\cos t_1 \gD x - \xi'(t_1)\gD \xi}l \gD\xi - l\xi'(t_1) = \frac{\partial}{\partial t_1}(l\Delta \xi).
\endalign$$
We omit the proof of the second equality since it is very similar.
\qed
\enddemo
\proclaim{Corollary \Tlabel\EisCinfty} $E_1$ is $C^\infty$ on 
$[-\frac\pi2,\frac\pi2]^2\bs\set{(-\frac\pi2,\frac\pi2),(\frac\pi2,-\frac\pi2),(0,0)}$
\endproclaim
\demo{Proof} Fix $(\ga_0,\gb_0)\in [-\frac\pi2,\frac\pi2]^2\bs\set{(-\frac\pi2,\frac\pi2),(\frac\pi2,-\frac\pi2),(0,0)}$
and let $N$ and $E^*$ be as in the proof above.  Then $E^*$ is $C^\infty$ on $Q(N)$, an open neighborhood of $(\ga_0,\gb_0)$.
The desired conclusion is now a consequence of the Claim in the above proof.
\qed\enddemo
%
\Section {Proof of Condition (\horseman)}
In this section, we prove condition (\horseman); namely that for every
$\ga\in[-\frac\pi2,\frac\pi2]$ there exists $\gb_\ga^*$,
with $\abs{\gb_\ga^*}\leq \frac\pi2-\Psi$, such that
$$
\text{sign}\pr{\frac{\partial}{\partial\gb}E_1(\ga,\gb)}=\text{sign}(\gb - \gb_\ga^*)\text{ for all }
\gb\text{ satisfying }\abs\gb\leq\frac\pi2\text{ and }\abs{\gb-\ga}<\pi.
\tag\Elabel\repeathorseman$$
With Theorem \inject\ in view, we treat the mapping $Q$ as a bijection
between $U$ and $Q(U)$, which (by Corollary \quickrain) contains $[-\frac\pi2,\frac\pi2]^2\bs\set{(0,0)}$.
Let $\ga\in[-\frac\pi2,\frac\pi2]$ be fixed. 
For the sake of clarity our proof is broken into three $\ga$-dependent cases.

\noindent{\bf Case 1:} $0<\ga\leq \frac\pi2$.

Set $B=[-\frac\pi2,\frac\pi2]\bs\set{\ga-\pi}$.
It follows from Corollary \EisCinfty\ that the function $\gb\mapsto E_1(\ga,\gb)$ is $C^\infty$ on $B$,
and, from Theorem \Tconsistent, we have
that $\frac{\partial}{\partial \gb} E_1(\ga,\gb)=\frac12 \gk_b(c_1(\ga,\gb))$.
Note that if $(t_1,t_2)=Q^{-1}(\ga,\gb)$, then $R_{[t_1,t_2]}$ is directly similar to $c_1(\ga,\gb)$,
and consequently 
$\text{sign}\pr{\frac{\partial}{\partial\gb}E_1(\ga,\gb)}=\text{sign}(\sin t_2)$ since the signed
curvature of $R(t)$ is $\gk(t)=2\sin t$.
\newline
\epsfysize=4.3truecm \epsffile{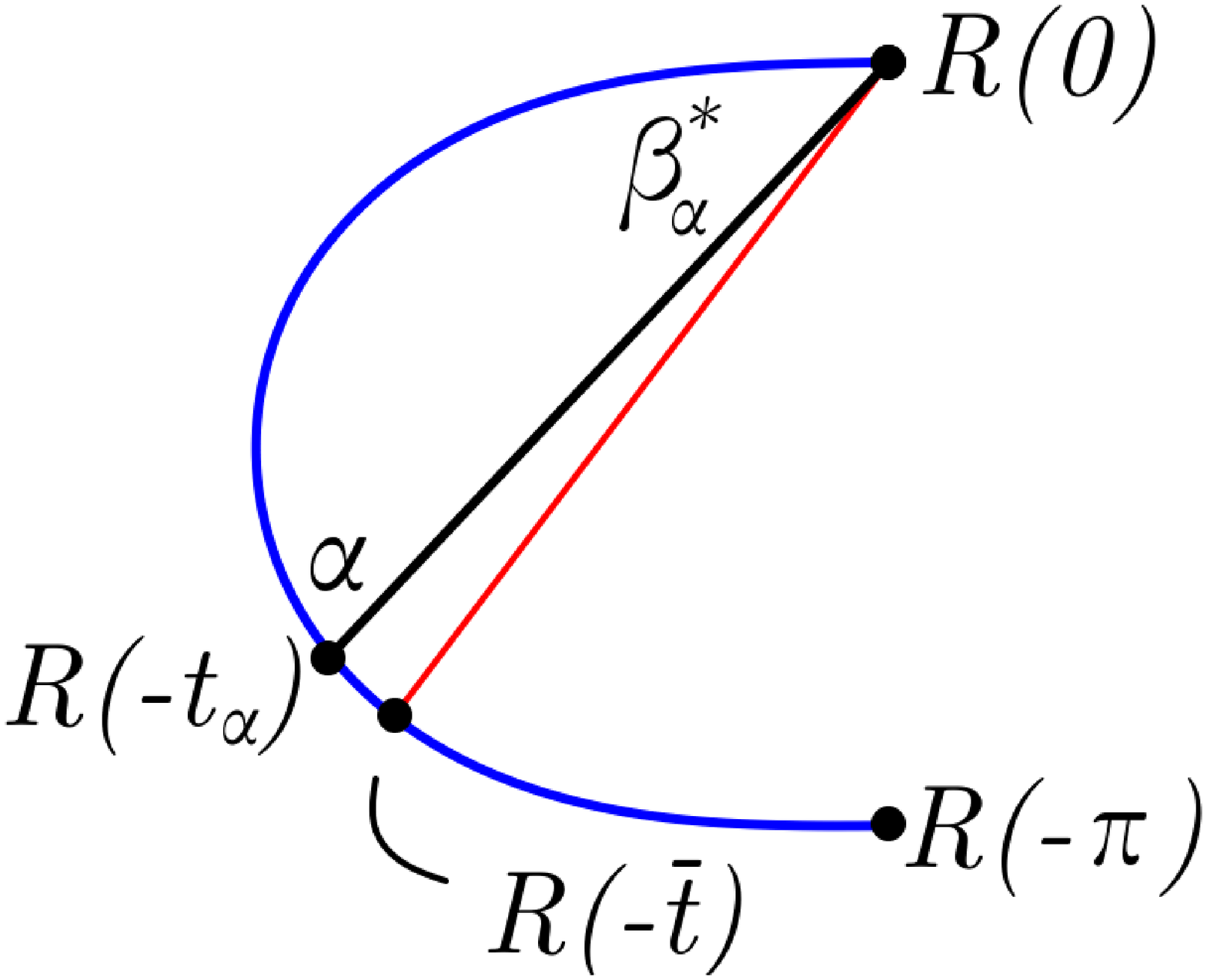}
\hskip2truecm  \epsfysize=4.3truecm \epsffile{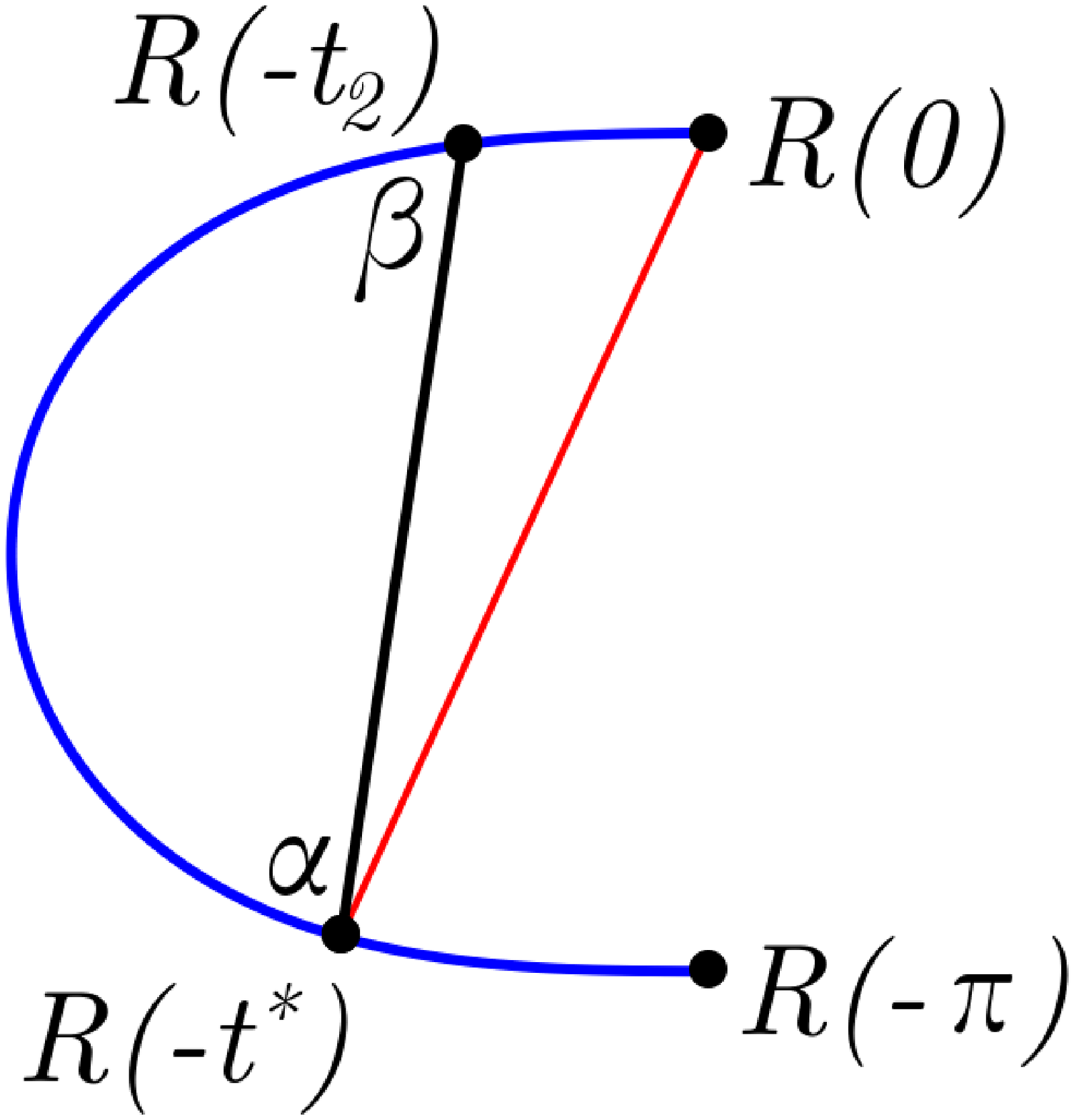}
\newline
{\bf Fig.\ \Flabel\figfourteen} \ the parameter $-t_\ga$ \hskip2truecm 
{\bf Fig.\ \Flabel\figfourteenb} \ the parameter $-t_2$\newline
By Theorem \detq\ and symmetry, there exists a unique $-t_\ga\in[-\tbar,0)$ such that $\ga(-t_\ga,0)=\ga$. Set
$\gb_\ga^*:=\gb(-t_\ga,0)<0$ and note that $R_{[-t_\ga,0]}$ (see Fig.\ \figfourteen) has chord angles $(\ga,\gb^*_\ga)$ while
$\text{sign}\pr{\frac{\partial E_1}{\partial\gb}(\ga,\gb_\ga^*)}=\text{sign}(\sin 0)=0$. Furthermore, we have
$\abs{\gb_\ga^*}=\abs{\ga(0,t_\ga)}\leq\abs{\ga(0,\tbar)}=\frac\pi2-\Psi$, and
it is shown in \cite{\BJone, Lemma 6.3} that
$\abs{\gb_\ga^*}=\abs{\ga(0,t_\ga)}<\gb(0,t_\ga)=\ga$. \newline
{\bf Claim:} If $\gb\in B$ is such that $\frac{\partial E_1}{\partial\gb}(\ga,\gb)=0$, then $\gb=\gb_\ga^*$.
\demo{proof} 
Assume $\gb\in B$ is such that $\frac{\partial E_1}{\partial\gb}(\ga,\gb)=0$. Set $(t_1,t_2)=Q^{-1}(\ga,\gb)$. Then
$t_2$ equals either $0$ or $\pi$ (since $\sin t_2=0$ and $(t_1,t_2)\in U$). If $t_2=0$, then $(t_1,t_2)\in U_0$ 
and it follows from
Theorem \detq\ and symmetry that $t_1=-t_\ga$ and hence $\gb=\gb_\ga^*$.  On the other hand, if $t_2=\pi$ then 
$(t_1,t_2)\in U_0$ and it follows that $\ga=\ga(t_1,t_2)<0$ which is a contradiction; hence the claim.
\enddemo
Note that $R_{[-t_\ga,t_\ga]}$ has chord angles
$(\ga,\ga)$ and hence
$\text{sign}\pr{\frac{\partial E_1}{\partial\gb}(\ga,\ga)}=\text{sign}(\sin t_\ga)>0$.  Since $\ga>0>\gb_\ga^*$, 
it follows from continuity
that $\text{sign}\pr{\frac{\partial E_1}{\partial\gb}(\ga,\gb)}>0$ for $\gb\in B$ with $\gb>\gb_\ga^*$.\newline
Now, in order to complete the proof (of Case I), it suffices to show that there exists $\gb\in B$ 
such that $\text{sign}\pr{\frac{\partial E_1}{\partial\gb}(\ga,\gb)}<0$. 
Since $\ga(-t^*,0)=\gb(0,t^*)>\frac\pi2\geq\ga$, it follows that there exists $-t_2\in(-t^*,0)$ such that
$\ga(-t^*,-t_2)=\ga$. Set $\gb:=\gb(-t^*,-t_2)<0$ (see Fig.\ \figfourteenb). It is easy to verify that $\abs\gb<\frac\pi2$
and therefore $\gb\in B$. Note that
$\text{sign}\pr{\frac{\partial E_1}{\partial\gb}(\ga,\gb)}=\text{sign}(\sin (-t_2)<0$. This completes
the proof for Case I.\newline

\noindent{\bf Case II:} $-\frac\pi2\leq\ga<0$
This case follows from Case I and the symmetry $E_1(\ga,\gb)=E_1(-\ga,-\gb)$.

\noindent{\bf Case III:} $\ga=0$.

Set $\gb_0^*:=0$. It is shown in Theorem \Tconsistent\ that $\frac{\partial E_1}{\partial\gb}(0,0)=0$.\newline
{\bf Claim:} If $\gb\in [-\frac\pi2,\frac\pi2]$ is such that 
$\frac{\partial E_1}{\partial\gb}(0,\gb)=0$, then $\gb=0$.
\demo{proof} 
By way of contradiction, assume $\gb\in [-\frac\pi2,\frac\pi2]\bs\set{0}$ is such that
$\frac{\partial E_1}{\partial\gb}(0,\gb)=0$. Set $(t_1,t_2)=Q^{-1}(0,\gb)$. Then
$t_2$ equals either $0$ or $\pi$. If $t_2=0$, then $t_1\in[-\pi,0)$ and it follows that $\ga>0$, which is
a contradiction. On the other hand, if $t_2=\pi$, then  $t_1\in[0,\pi)$ and it follows that $\ga<0$, which 
again is a contradiction; hence the claim.
\enddemo
The symmetry $E_1(0,\gb)=E_1(0,-\gb)$ ensures that 
$\frac{\partial E_1}{\partial\gb}(0,-\gb)=-\frac{\partial E_1}{\partial\gb}(0,\gb)$ and therefore it suffices
to show that $\frac{\partial E_1}{\partial\gb}(0,\gb)>0$ for all $\gb\in (0,\frac\pi2]$.
Define $g(\gb):=E_1(0,\gb)$, $\gb\in [0,\frac\pi2]$ so that 
$g'(\gb)=\frac{\partial E_1}{\partial\gb}(0,\gb)$.
Then $g$ is continuous on $[0,\frac\pi2]$
and is $C^\infty$ on $(0,\frac\pi2]$. 
It follows from the claim that $\text{sign}\pr{g'}$ is nonzero and constant on  $(0,\frac\pi2]$.
If $\text{sign}\pr{g'}=-1$ on 
$(0,\frac\pi2]$, then we would have $E_1(0,\frac\pi2)<E_1(0,0)=0$, which is a contradiction; therefore $\text{sign}\pr{g'}=1$ on 
$(0,\frac\pi2]$ and this completes the proof of the final case.
%
\Section {Appendix }
%
The goal of this section is to prove Theorem \inject. The proof that
$Q$ is injective on $U$ uses ideas from the proof of the
Hadamard-Caccioppoli theorem, which states [{\AmProdi}, Th. 1.8,
page 47]

\proclaim {Theorem (Hadamard-Caccioppoli)} Let $M,N$ be metric
spaces and $F\in C(M,N)$ be proper and locally invertible on all
of $M$. Suppose that $M$ is arcwise connected and $N$ is simply
connected. Then $F$ is a homeomorphism from $M$ to $N$.
\endproclaim

 Unfortunately, not all the
conditions of the Hadamard-Caccioppoli theorem are satisfied, for
$Q$ is not locally invertible on $\gamma _0$.  To remedy this we
are going to use results needed in the proof the
Hadamard-Caccioppoli theorem [{\AmProdi}, Th. 1.6, page 47].

Let $M,N$ be metric spaces. For a map $F\in C(M,N)$ denote by
$\Sigma=\{u\in M: F \hbox{\ is not locally invertible at\ }u\}$
the singular set of $F$ and for $v\in N$ denote by $[v]$ the
cardinal number of the set $F^{-1}(v)$.

\proclaim {Theorem \Tlabel \pproper} [{\AmProdi}, Th. 1.6, page 47]
Let  $F\in C(M,N)$
be proper.  Then $[v]$ is  constant on every connected component
of $N-F(\Sigma)$.
\endproclaim

In our case $M=U\cup \gamma _0$, $F=Q$ and $N=Q(M)$. The
properties of $Q$ are summarized in Proposition \injectgamma.

Let us recall from Section 6 that $M$ is topologically an annulus,
the boundary $\partial M$ consists of two curves $\gamma _0$ and
the staircase curve $\gamma$ depicted in Fig.\ \figten\ (b). Since $M$ is
compact $Q$ will be proper ($Q^{-1}(K)$ is compact if $K$ is
compact). $Q(\gamma _0)=\{(0,0)\}$ and $Q(\gamma)=\Gamma $
(depicted in Fig.\ \figten\ (a)) is a simple closed curve ((iv) of
Proposition \injectgamma) and $Q$ is injective on $\gamma$.

First we will show that $Q$ maps $M$ onto the union of the
interior of $\Gamma$ and $\Gamma$ and it maps the interior of $M$
onto the interior of $\Gamma$ minus the point $(0,0)$.

Since $\Gamma$ is a Jordan curve it has an interior. Let us denote
by $N_0=\{\hbox{ interior of\ } \Gamma\}\cup \Gamma$.
\proclaim {Proposition \Tlabel \intgamma} $N=N_0$ and
$Q(M-\{\gamma, \gamma _0\})=N-\big(\Gamma \cup\{(0,0)\}\big).$
\endproclaim

\demo {Proof} Claim 1. If $x\in \hbox {int}M$, then $Q(x)\in N_0$
(where int$M=M-\{\gamma,\gamma_0\}$). Suppose it is not true and
there is a point $x\in \hbox {int}M$ such that $Q(x)\notin N_0$.

By Proposition \injectgamma\ the Jacobian of $Q$  is not zero at the points
of $\hbox{int}M$ therefore it is an open mapping, that is
$Q(\hbox{int}M)$ is an open set. The indirect assumption means
that $Q(\hbox{int}M)$ has a point outside $N_0$ and since it is a
bounded set (obvious from the definition of $Q$) it must have a
boundary point $y$ outside $N_0$. Let $x_n\in \hbox {int}M$ be a
sequence such that $\lim Q(x_n)=y$. Passing on to a subsequence if
necessary we can assume that $x_n$ is convergent with $\lim
x_n=x^*$. Clearly $x^*\notin \{\gamma, \gamma _0\}$ since
$Q(\gamma _0)=(0,0)$ and $Q(\gamma )=\Gamma $. Since the Jacobian
of $Q$ at $x^*$ is not zero therefore it maps an open neighborhood
of $x^*$ onto an open neighborhood of $Q(x^*)=y$, which is a
contradiction.

Claim 2. If $x\in \hbox {int}M$, then $Q(x)\in N_0-\{\Gamma,
(0,0)\}$. Item (iii) of Proposition \injectgamma\ states that $Q(x)\ne (0,0)$.
Since $Q$ is a local diffeomorphism at $x$ ((ii) Proposition \injectgamma) if
$Q(x)\in \Gamma $ that would imply that there is a point $y\in
\hbox {int}M$ near $x$ such that $Q(y)\notin N_0$. That would be
in contradiction with Claim 1.

Claim 3. The map $Q:M\to N_0$ is onto. Since $Q(\gamma _0)=(0,0)$,
$Q(\gamma)=\Gamma$ and $Q(M)$ compact if $Q(M)\ne N_0$ there has
to be a boundary point $u$ of $Q(M)$ such that $u\in Q(M)-(\Gamma
\cup \{(0,0)\})$. To find $u$ one has to connect a point inside
the image different from $(0,0)$ to a point of $N_0$ which is
outside $Q(M)$ with a curve avoiding both the point $(0,0)$ and
the curve $\Gamma$. On this curve one can find $u$.

 Let $x\in \hbox {int}M$ be a point such that
$Q(x)=u$. Since $Q$ is a local diffeomorphism at $x$, $u=Q(x)$
cannot be a boundary point of the image. This leads to a
contradiction and the claim is proved.

Therefore we have $N=N_0$ and  the Proposition is proved. \qed
\enddemo
\proclaim {Proposition \Tlabel \injgamma} For any $v\in
\Gamma$, $Q^{-1}(v)$ consists of one single point.
\endproclaim
\demo {Proof} Since $Q(\gamma)=\Gamma$ it is clear that
$Q^{-1}(v)$ is not empty. From the Previous Proposition it follows
that if $x,y\in Q^{-1}(v)$, then both $x,y\in \gamma$. Since $Q$
is injective on $\gamma$ ((iv) Theorem 7.8) it implies that
$x=y$\qed \enddemo
To show that the singular set of $Q$ (the set of points where $Q$
is not locally invertible) is $\gamma_0$ only, we need the
following:

\proclaim {Proposition \Tlabel \localgamma} $Q$ is
locally invertible at every point of $\gamma$.
\endproclaim

The proof will rely on the fact [{\chwuho}, Lemma 3, page 239]
that proper local homeomorphisms are covering maps, therefore they
have the unique path-lifting property. The precise statement is as
follows.

\proclaim {Proposition \Tlabel\proper} [{\chwuho}, Lemma
3] Let $X,Y$ be two Hausdorff spaces and let $Y$ be pathwise
connected. Any surjective, proper local homeomorphism $f:X\to Y$
must be a covering projection.
\endproclaim
\demo {Proof of Proposition {\localgamma}} Let $x\in \gamma$ be
any point and set $y=Q(x)\in \Gamma$. Choose $\epsilon <
dist(\Gamma , (0,0))$ small enough such that $N\cap B(y,\epsilon)$
is connected therefore simply connected. Here $B(y,\epsilon)$
denotes the closed ball of radius $\epsilon$ centered around $y$.
Since the boundary of $N$ is a piecewise differentiable curve one
can find such $\epsilon$.

Let $\delta >0$  be chosen such that $Q(M\cap B(x,\delta)\subset
B(y,\epsilon)$. Such $\delta$ exists because $Q$ is continuous at
$x$. It is enough to show that $Q$ is one-to-one on $M\cap
B(x,\delta)$ since a continuous, one-to-one map between compact
subsets of $\Bbb R^2$ has a continuous inverse.

Suppose it is not true. Then there are points $x_1\ne x_2\in M\cap
B(x,\delta)$ such that $Q(x_1)=Q(x_2)$. From the previous
propositions (Proposition {\injgamma} and Proposition {\intgamma})
it follows that $x_1,x_2\notin \gamma $. Let $h:[0,1]\to \hbox
{int}M\cap B(x,\delta)$ be a curve connecting $x_1$ to $x_2$ and
set $g=Q(h)$. Then $g$ is a closed curve in $\hbox {int}N\cap
B(y,\epsilon)$ which is simply connected. Note that $\hbox
{int}N\cap B(y,\epsilon)$ does not contain $(0,0)$. Therefore,
there is a homotopy $H:[0,1]\times[0,1]\to N\cap B(y,\epsilon)$
such that $H(0,t)=g(t)$ and $H(s,0)=H(s,1)=H(1, t)=Q(x_1)$ for all
$s,t\in [0,1]$.

The map $Q:\hbox {int}M\to \hbox {int}N-\{(0,0)\}$ is proper
(Proposition {\intgamma}) and since the Jacobian of $Q$ is not
zero, it is a local homeomorphism, hence by Proposition {\proper}
it is a covering map. This means that we can lift $H$ (the image
of $H$ avoids the point $(0,0)$) to a homotopy $\bar
H:[0,1]\times[0,1]\to \hbox {int}M$ with the property that $Q(\bar
H(s,t))=H(s,t)$. This implies that $\bar H(s,0)=x_1$ and $\bar
H(s,1)=x_2$ for all $s\in [0,1]$, therefore we have curve $t\to
\bar H(1,t)$ in $\hbox {int}M$ connecting $x_1$ to $x_2$ such that
the image of this curve by $Q$ is one point $Q(x_1)$. This
contradicts the fact that $Q$ is a local homeomorphism on $\hbox
{int}M$. \qed \enddemo
 \demo {Proof of Theorem \inject} We have already proved that  $Q$ is proper.
Since $Q$ is locally invertible at the points of $\hbox {int}M$
(the Jacobian of $Q$ is not zero) and at the points of $\gamma$
(Proposition {\localgamma}) we see that the singular set of $Q$ is
$\gamma _0$ only and since $Q(\gamma _0)=\{(0,0)\}$ from Theorem
{\pproper} we obtain that for all $v\in N-\{(0,0)\}$ $[v]$ is
constant. Since $Q$ is injective on $\gamma$ and if $u\in
M-\{\gamma, \gamma _0\}$ then $Q(u)\notin \Gamma$ we obtain that
$[v]=1$ for all $v\in \Gamma$, therefore $[v]=1$ for all $v\in
N-\{(0,0)\}$. This means that $Q$ is injective on $U$ and the
proof of Theorem \inject\ is complete. \qed
\enddemo

%

\enddocument

%% file: SecThEq.tex
\newcount\sectionno\sectionno=0
\newcount\equationno\equationno=0
\newcount\theoremno\theoremno=0
\newcount\figureno\figureno=0

\font\bigfont=cmbx10
\def\Heading#1{\bigskip\bigskip\goodbreak\centerline{\bigfont #1}\nobreak\bigskip\nobreak}
\def\Section#1{\global\advance\sectionno by 1\relax%
   \equationno=0\theoremno=0\Heading{\the\sectionno. #1}}

\newwrite\aux
\immediate\openout\aux=\jobname.aux

\def\formeq{\the\sectionno.\the\equationno}
\def\formth{\the\sectionno.\the\theoremno}
\def\formfig{\the\figureno}

\def\Elabel#1{\global\advance\equationno by 1 \formeq%
\immediate\write\aux{\def\string#1{\formeq}}%
\global\edef#1{\formeq}}

\def\Tlabel#1{\global\advance\theoremno by 1 \formth%
\immediate\write\aux{\def\string#1{\formth}}%
\global\edef#1{\formth}}

\def\Flabel#1{\global\advance\figureno by 1 \formfig%
\immediate\write\aux{\def\string#1{\formfig}}%
\global\edef#1{\formfig}}